\documentclass[12pt]{article}
\usepackage[utf8]{inputenc}
\usepackage[margin=1in]{geometry}
\usepackage{enumerate}
\usepackage{graphicx,float}
\usepackage{amssymb,MnSymbol,dsfont}
\usepackage{doi}
\usepackage{color}
\usepackage{hyperref}
\usepackage{cleveref}
\usepackage{algorithm}
\usepackage{algorithmic}
\usepackage{subcaption}

\newcommand{\RR}{\mathds{R}}
\newcommand{\TT}{\mathds{T}}

\newtheorem{thm}{Theorem}
\newtheorem{lem}{Lemma}
\newtheorem{asp}{Assumption}
\newtheorem{defi}{Definition}
\newtheorem{rmk}{Remark}
\DeclareMathOperator{\ad}{ad}
\newenvironment{pf}{\noindent\textbf{Proof.}}{\hfill$\square$}

\newcommand\numberthis{\addtocounter{equation}{1}\tag{\theequation}}
\DeclareUnicodeCharacter{2212}{\ensuremath{-}}

\title{Exponential time differencing for matrix-valued dynamical systems\thanks{This work was supported by the Engineering and Physical Sciences Research Council (EPSRC) through the Mathematics of Systems II Centre for Doctoral Training at the University of Warwick (reference EP/S022244/1). T.G. acknowledges the support received from the EPSRC projects EP/T011866/1 and EP/V013319/1. For the purpose of open access, the author has applied a Creative Commons Attribution (CC BY) licence to any Author Accepted Manuscript version arising from this submission.}}

\date{\today}

\author{Nayef Shkeir\\
\small Mathematics Institute, University of Warwick, Coventry, CV4 7AL, United Kingdom\\
\and Tobias Grafke\\
\small Mathematics Institute, University of Warwick, Coventry, CV4 7AL, United Kingdom}

\begin{document}

\maketitle

\begin{abstract}
    Matrix evolution equations occur in many applications, such as dynamical Lyapunov/Sylvester systems or Riccati equations in optimization and stochastic control, machine learning or data assimilation. In many such problems, the dominant stability restriction is imposed by a stiff linear term, making standard explicit integrators impractical. Exponential time differencing (ETD) is known to produce highly stable numerical schemes by treating the linear term in an exact fashion. In particular, for stiff problems, ETD methods are the methods of choice. We extend ETD to matrix-valued evolution equations of the form $\dot Q = LQ + QR + N(Q,t)$ by deriving explicit matrix-ETD (METD) schemes. When $L$ and $R$ commute, we construct an explicit $p$-th order METD$p$ family and prove order-$p$ global convergence under standard assumptions; for the non-commuting case, we develop a Baker–Campbell–Hausdorff (BCH)-based extension. This allows us to produce highly efficient and stable integration schemes. We demonstrate efficiency and applicability on stiff PDE-derived and large-scale matrix dynamics, including an Allen–Cahn system, turbulent jet fluctuation statistics, and continuous graph neural networks. We further show that the scheme is more accurate, stable, and efficient than competing schemes in large-scale high-rank stiff systems.
\end{abstract}
\section{Introduction}

Matrix-valued dynamical systems play a crucial role in understanding complex phenomena in the sciences. They often occur as differential Lyapunov, differential Sylvester or matrix-Riccati equations in the context of optimization, but are also present in stochastic control and data assimilation. Differential Lyapunov equations appear in stability analysis for continuous linear time-varying systems \cite{amato2014finite} where analysis and controllability depend on the solutions of the Lyapunov equations. Differential Lyapunov/Sylvester equations also have important applications in optimal control \cite{athans2013optimal}, model reduction \cite{gugercin2008h_2} and machine learning \cite{pmlr-v119-xhonneux20a}. Matrix-Riccati equations play a significant role in the context of continuous-time Kalman filters, which are widely utilized in numerous scientific and engineering disciplines for state estimation and control of dynamical systems \cite{kalman1961new,ricatti_book}.

Matrix-valued dynamical equations usually come with tight stability constraints, for example, because of the presence of high-order linear differential operators (e.g.~diffusion or (hyper-)viscosity), or poor conditioning in the (quasi-)lin\-ear dynamics \cite{cho1996emergence,danilov2000quasi}. Moreover, many real-world systems have multiple time scales or are oscillatory, making the system stiff~\cite{hingham_stiff}. For systems with scale separation, semi-Lagrangian methods have been popular in physics applications \cite{PRUSA20081193}. 

There has been some effort to solve matrix-valued dynamical systems numerically, such as backward differentiation formula (BDF) methods \cite{benner2004bdf}, Parareal based algorithms \cite{parareal_ricatti} and Krylov subspace methods \cite{behr_kyrlov,kyrlov_2}. There has also been some preliminary work on the use of exponential integrators for Matrix-Riccati equations \cite{LI2021113360}. However, many existing approaches become infeasible for large and stiff matrix-valued systems. Recently, projected exponential methods  have been developed that apply Krylov subspace techniques to approximate $\varphi$-functions of the Sylvester operator directly, without requiring commutativity of the coefficient matrices \cite{carrel2023projected}. The present work takes a complementary approach: we exploit commutativity (common in stochastic applications) to derive explicit closed-form schemes where $\varphi$-functions are evaluated on the coefficient matrices directly, rather than on the full Sylvester operator.

For scalar systems, exponential time differencing (ETD) methods provide accurate numerical solutions to stiff dynamical problems \cite{cox-matthews:2002, friedli1978generalized} and have been applied to a variety of problems involving dynamical systems modeled by ordinary and partial differential evolution equations \cite{etd_example_1,ZonostrophicInstability,turb,frisch2008hyperviscosity,ju2015fast}. These schemes have been shown to be especially suited to quasi-linear systems which can be split into a linear part that contains much of the stiffness of the dynamics and a nonlinear part \cite{Minchev2005ARO,ZonostrophicInstability}. 

Given these two developments, it is natural to ask if ETD methods can be extended to matrix-valued evolution equations by making use of the particular structure present in such problems. The main objective of this present work is to present a method to generalize ETD schemes to matrix-valued problems. Similar to the case of standard scalar ETD problems, it turns out that the key idea of the approach, namely to integrate the linear part of the problem exactly, can be carried over to matrix-valued problems. The remaining nonlinear integral term, however, must be treated with more care than in the scalar case, since non-commutativity gives rise to complications, leading to additional terms in the numerical scheme. 

The structure of this paper is as follows: We first review the basic idea of ETD algorithms in their extension to matrix-valued equations. A particular focus is on commutator relations that make the study of ETD algorithms in the matrix case fundamentally different from the scalar case. We then introduce a first-order matrix-ETD (METD) scheme, and two second-order algorithms, a METD multistep scheme and a METD Runge-Kutta scheme. We also introduce a general METD$p$ multistep construction that achieves order-$p$ accuracy by systematically retaining nested commutator corrections, and provide a complete error analysis establishing both local consistency and global convergence of order $p$ (Theorems~\ref{thm:metdp_local} and~\ref{thm:metdp_global}). We then extend the methodology to non-commuting coefficient matrices using the Baker--Campbell--Hausdorff series. After the theoretical analysis, we illustrate the methods on (i) a low-dimensional differential Lyapunov test that isolates commutator effects, (ii) a stiff nonlinear Allen–Cahn benchmark arising from PDE discretization, and (iii) a large-scale Lyapunov system governing turbulent jet fluctuations. Finally, we discuss an application to continuous graph neural networks involving non-square matrices and non-commuting coefficients, demonstrating how the extensions can be combined. All schemes derived in this paper are explicit.

\section{Matrix Exponential Time Differencing}

Let us consider a matrix evolution equation for $Q(t) \in \mathbb{R}^{n \times n}$ given by 
\begin{equation}
    \label{eq:ODE}
    \dot Q = LQ + QR + N(Q,t),\quad Q(0) = Q_0\,,
\end{equation}
with linear operators $L, R \in \mathbb{R}^{n \times n}$ multiplied once from the left and from the right, respectively, against $Q$, and with $N(Q,t): \mathbb{R}^{n\times n}\times [0,T]\to  \mathbb{R}^{n \times n}$ being potentially nonlinear. Note that, in what follows, we define the commutator $[A,B] := AB - BA\,$ for $A,B \in \mathbb{R}^{n \times n}$ as usual. It is also important to note that the linear operators ($L$ and $R$) are in the same space as the unknown $Q$.

The first step toward a solution of the problem (\ref{eq:ODE}) is to solve the linear problem. In the scalar case, this is obtained by using an integrating factor. In the matrix case at hand, we can multiply $Q$ from the left by $e^{-tL}$ and right by $e^{-tR}$ and differentiate 
to obtain
\begin{align*}
  \frac{d}{dt}\left(e^{-tL} Q e^{-tR} \right) 
  &= e^{-tL} \dot Q e^{-tR} - e^{-tL}  LQ e^{-tR} - e^{-tL} QR e^{-tR}\\
  &= e^{-tL} N(Q,t) e^{-tR} \,, &&
\end{align*}

an equation we can directly integrate. Here $e^{A}$ for $A \in \RR^{n\times n}$ is the matrix exponential. The solution of (\ref{eq:ODE}) in the time interval $[t_n,t_{n+1}]$ where $t_{n+1} = t_n + \Delta t$ and $t_0=0$, is given by the variation-of-constants formula  
\begin{equation}
  \label{eq:etd-int}
  Q(t_{n+1}) = e^{\Delta t L} Q(t_n)\, e^{\Delta t R} \, \nonumber + \int_0^{\Delta t} e^{(\Delta t-\tau) L} N(t_n + \tau) e^{(\Delta t-\tau) R}\,d\tau\,, \numberthis
\end{equation}

where  we need to evaluate the integral in (\ref{eq:etd-int}), which, if the commutator $[N,e^{sR}]$ is non-zero, is not trivial (equally if $[N,e^{sL}] \neq 0$). With classic ETD methods, we try to approximate the nonlinear term in the integral using an algebraic polynomial. For now, we consider the case $[L,R] = 0$, and will drop this simplification only in section~\ref{sec:extensions_non_comm}. Note that many (but not all) problems of practical relevance satisfy this condition (a common case is, e.g.,~$R=L^T$ and $L$ normal). Throughout this paper, we adopt the notation $N(Q(t),t) = N(t)$ and $N(t_n)=N_n$ for simplicity.

We also define $\varphi$ functions commonly seen in the exponential integrators literature \cite{LUAN2014168, hochbruck2010exponential} to simplify notation. Let $\theta = \tau / \Delta t \in [0,1]$ and we then define these functions as 
\begin{equation}
    \label{eq:varphi_def}
    \varphi_k(A) =\int_0^{1} e^{(1 - \theta)A} \frac{\theta^{k-1}}{(k-1)!}\, d\theta,  \qquad  k\geq 1 \, ,
\end{equation}
with $\varphi_0(A) = e^{A}$ for some matrix $A \in \RR^{n\times n}$. We also note that 
\begin{equation}
    \varphi_{k+1}(A) = A^{-1}[\varphi_k(A) - \varphi_k(0)], \qquad \varphi_k(0) = \frac{\mathds 1}{k!}\,,
\end{equation}
where $\mathds 1_n$ is the identity on $\RR^{n\times n}$ and we drop the subscript when unambiguous from context. For most of the numerical examples we present in this paper, we can compute our $\varphi$ functions directly and in cases where $A$ is ill-conditioned, we introduce regularization. There is a breadth of literature on computing $\varphi$ functions efficiently such as Krylov subspace methods \cite{Niesen_2012} and the Padé method with the scaling and squaring strategy \cite{AlMohy2009970}.

\begin{rmk}
The left- and right-multiplication operators $\mathcal{L}(X) = LX$ and $\mathcal{R}(X) = XR$ always commute as linear operators, regardless of whether $LR = RL$. The assumption $[L,R] = 0$ is imposed not for the validity of the fundamental solution $e^{tL}Q_0 e^{tR}$, but to enable the simplification $e^{sL}e^{sR} = e^{s(L+R)}$. This allows coefficients to be expressed compactly in terms of $\varphi_k$.
\end{rmk}

\subsection{Commutator series}
We can multiply the commutator $[N,e^{sR}] $ from the left by $e^{sL}$
\begin{equation*}
    e^{sL}[N,e^{sR}] = e^{sL}Ne^{sR} - e^{s(L+R)}N\,,
\end{equation*}
where $e^L e^R = e^{L+R}$ only because we demand commutativity. This means we can rewrite our integral in (\ref{eq:etd-int}) as
\begin{align}
 \label{eq:integ}
    \int_0^{\Delta t} e^{(\Delta t-\tau)L} N(t_n+\tau) e^{(\Delta t-\tau)R}\,d\tau  \nonumber =  &\int_0^{\Delta t} e^{(\Delta t-\tau)(L+R)} N(t_n+\tau)\,d\tau \nonumber  \,  \\ &+ \int_0^{\Delta t} e^{(\Delta t-\tau)L} [N(t_n+\tau),e^{(\Delta t-\tau)R}]\,d\tau.
\end{align}
We can derive numerical schemes via a truncation of the commutator by considering the power expansion
\begin{equation*}
  e^{sR} = \sum_{k=0}^\infty \frac{s^k}{k!}R^k\,,
\end{equation*}
which gives us the commutator series expansion 
\begin{equation*}
  \label{eq:commutator-series}
  [N, e^{sR}] = \sum_{k=0}^\infty \frac{s^k}{k!} [N,R^k] = s[N,R] + \tfrac12 s^2 [N,R^2] + \dots \, .
\end{equation*}
Our goal is to approximate the integral up to a certain order and considering that $s \in [0,\Delta t]$, truncating the commutator series at a certain order of $s$ will allow us to 
carry out this expansion.

\subsection{First-order Matrix Exponential Time Differencing Scheme, METD1}

With all the building blocks in place, we can quickly arrive at the first-order scheme, using
{
\begin{align*}
    \label{eq:comm_integ}
    \int_0^{\Delta t} e^{(\Delta t-\tau)L}  N(t_n+\tau) e^{(\Delta t-\tau)R}\,d\tau  =  & \, \Delta t\int_0^{1} e^{(1-\theta)\Delta t(L+R)} N(t_n+\Delta t \theta)\,d\theta  \,  \\ &+ \Delta t^2 \int_0^{1} (1-\theta) e^{(1-\theta)\Delta t L} [N(t_n+\Delta t \theta),R]\,d\theta    \\ &+ \mathcal O(\Delta t^3) \, \numberthis
   \\ \approx &\, \Delta t\int_0^{1} e^{(1-\theta)\Delta t(L+R)} \,d\theta N(t_n) 
   \\ = \, &(L+R)^{-1}\left(e^{\Delta t (L+R)} - \mathds{1}\right) N(t_n) \,,
\end{align*}
}
where we approximate $N$ to first-order as a constant between $t_{n}$ and $t_{n+1}$ (in the following, index $n$ represents the temporal index for discrete time steps). Details of the computation of the integral are given in \ref{ETD1_deriv}. The commutator integral in (\ref{eq:comm_integ}) is $\mathcal O(\Delta t^2)$ which justifies dropping the commutator integral for the METD1 scheme. 

Summarizing these calculations, we have obtained a first-order METD scheme as
\begin{equation}
\label{eq:etd1-eqn}
  Q_{n+1} = \varphi_0(\Delta t L)\, Q_{n}\,  \varphi_0(\Delta t R) + \Delta t \, \varphi_1(\Delta t(L+R))\, N(Q_{n}, t_n)\,,
\end{equation}
where by solving the corresponding integrals in (\ref{eq:varphi_def}), we have
\begin{align*}
  &\varphi_0(\Delta tA) = e^{\Delta t A}\; \text{,}\\ &\varphi_1(\Delta tA) = (\Delta t A)^{-1}\left[e^{\Delta t A} - \mathds{1}\right]\,.
\end{align*}

\subsection{Second-order Multistep scheme, METD2}
\label{sec:metd2}
In order to derive a second-order scheme, we proceed similarly to the scalar case \cite{cox-matthews:2002} and now approximate the nonlinear term as 
\begin{equation*}
    N(t_n + \tau) = N_{n} +\tau (N_{n}- N_{n-1})/\Delta t+ \mathcal O(\Delta t^{2}),
\end{equation*}
 which is a first-order polynomial in $\tau$, where $\tau$ is in the interval $0\leq \tau \leq \Delta t$. We carry out the same steps as for the first-order scheme, but now have to consider the commutator expansion in the integral in \eqref{eq:integ}. Defining the backward difference as $\nabla N_n := N_n - N_{n-1}$, we arrive at the following integrals 
\begin{align}
    \label{eq:ETD2_integ}
    \int_0^{\Delta t} e^{(\Delta t-\tau)L} N(t_n+\tau) e^{(\Delta t-\tau)R}\,d\tau =& \nonumber  \nonumber
    \Delta t\int_0^{1} e^{(1-\theta)\Delta t(L+R)} N(t_n+\Delta t \theta)\,d\theta \nonumber \\
    & +\Delta t^2 \int_0^{1} (1-\theta) e^{(1-\theta)\Delta t L} [N(t_n+\Delta t \theta),R]\,d\theta \nonumber \\ &+ \mathcal O(\Delta t^3) \\ 
    \approx & \, \Delta t \int_0^{1} e^{(1-\theta)\Delta t(L+R)} (N_n + \theta \nabla N_n)\,d\theta \,  \nonumber  \\
    & + \Delta t^2 \int_0^{1} (1-\theta) e^{(1-\theta)\Delta t L} [(N_n + \theta \nabla N_n),R]\,d\theta  \nonumber
\end{align}
to second-order. We treat these integrals in detail in \ref{app:METD_deriv}. Putting everything together, we obtain the METD2 scheme as
\begin{align}
\label{eq:etd2-eqn}
  Q_{n+1} = \,&\varphi_0(\Delta tL)\, Q_{n} \, \varphi_0(\Delta tR) + \Delta t \, \varphi_1(\Delta t(L+R))\, N_n \, + \nonumber \\&\Delta t \,\varphi_2(\Delta t(L+R))\, \nabla N_n + \Delta t^2 \,(\varphi_1(\Delta t L)- \varphi_2(\Delta t L))\, [N_n,R]\,, 
\end{align}
where by solving the corresponding integrals in (\ref{eq:varphi_def}), we have
\begin{align*}
    &\varphi_0(\Delta tA) = e^{\Delta t A}\; \text{,} \\ &\varphi_1(\Delta tA)  = (\Delta t A)^{-1}\left[e^{\Delta t A} - \mathds{1}\right]\; \text{,} \\ &\varphi_2(\Delta tA)  = (\Delta t A)^{-2}\left[e^{\Delta t A} - \Delta t A - \mathds{1}\right]\, .
\end{align*}

\subsection{Second-order Runge-Kutta scheme, METD2RK}

Multistep methods can be inconvenient to implement with only one initial value. Runge-Kutta (RK) methods avoid this problem, and further often have larger stability regions and smaller error constants compared to multistep methods. This comes at the price of a larger number of evaluations of the RHS per step. For a second-order RK method, we take an intermediate step given by 
\begin{equation*}
    A_{n} = \varphi_0(\Delta tL)\, Q_{n}\, \varphi_0(\Delta tR) + \Delta t \, \varphi_1(\Delta t(L+R))\, N(Q_{n}, t_n)\, ,
\end{equation*}
which is just the first-order scheme~(\ref{eq:etd1-eqn}). For the METD2RK scheme, we take the second-order approximation of the nonlinear
term and consider
\begin{equation*}
   N(t_n+\tau) = N(Q_{n}, t_n) +\tau (N(A_{n}, t_n + \Delta t) -  N(Q_{n}, t_n))/\Delta t+ \mathcal O(\Delta t^{2})\, ,
\end{equation*}
which, combined with \eqref{eq:ETD2_integ}, gives the METD2RK scheme as 
\begin{align*}
     Q_{n+1} =& \, A_{n} + \Delta t \, \varphi_2(\Delta t(L+R)) \big(N(A_{n}, t_n + \Delta t)- N(Q_{n}, t_n)\big) \\
     &+ \Delta t^2 \,(\varphi_1(\Delta t L)- \varphi_2(\Delta t L))\, [N(Q_{n}, t_n),R]\,.
\end{align*}

\section{A general METD\texorpdfstring{$p$}{p} construction}
\label{sec:metdp}

We now state a general fixed-step METD scheme of order $p$ for \eqref{eq:ODE} in the commuting case $[L,R]=0$.
The detailed derivation is given in Appendix~\ref{app:METDp_deriv}.
Throughout, we set $h=\Delta t$ and use the shorthand $N_n := N(Q_n,t_n)$.

Define the linear operator $\ad_R:\RR^{n\times n}\to\RR^{n\times n}$ by
\begin{equation}
    \ad_R(X):=[X,R]=XR-RX.
\end{equation}
The induced exponential satisfies the conjugation identity
\begin{equation}
\label{eq:ad_conj}
    e^{s\ad_R}(X)=e^{-sR}Xe^{sR}, \qquad s\in\RR,
\end{equation}
with a proof given in Appendix~\ref{app:METDp_deriv}.
Using \eqref{eq:ad_conj} and $[L,R]=0$, we can rewrite the integrand as
\begin{equation}
\label{eq:metdp_integrand_rewrite}
    e^{sL}Xe^{sR} = e^{s(L+R)}\,e^{s\ad_R}(X), \qquad s\ge 0.
\end{equation}
Hence, starting from the variation-of-constants formula \eqref{eq:etd-int}, we denote the nonlinear contribution over one step by
\begin{equation}
\label{eq:metdp_I_def}
    I := \int_0^h e^{(h-\tau)L}\,N(t_n+\tau)\,e^{(h-\tau)R}\,d\tau .
\end{equation}
Using \eqref{eq:metdp_integrand_rewrite} with $s=h-\tau$ yields the exact representation
\begin{equation}
\label{eq:metdp_integral_exact}
    I = \int_0^h e^{(h-\tau)(L+R)}\,e^{(h-\tau)\ad_R}\bigl(N(t_n+\tau)\bigr)\,d\tau.
\end{equation}
Let $\nabla$ denote the backward-difference operator acting on the sequence $(N_n)$ (just like in \ref{sec:metd2}),
\begin{equation*}
\nabla^0N_n=N_n,\qquad 
\nabla^m N_n=\sum_{\ell=0}^{m}(-1)^\ell\binom{m}{\ell}N_{n-\ell},\qquad m\ge 1.
\end{equation*}
As in classical ETD multistep schemes \cite{cox-matthews:2002}, we approximate $N(t_n+\tau)$ on $[0,h]$ by the degree-$(p-1)$ Newton backward interpolant,
\begin{equation}
\label{eq:metdp_interp}
    N(t_n+\tau)=\sum_{m=0}^{p-1}(-1)^m\binom{-\theta}{m}\,\nabla^m N_n + \mathcal O(h^p),
    \qquad \theta=\tau/h\in[0,1].
\end{equation}
Changing variables $\tau=h\theta$ in \eqref{eq:metdp_integral_exact} and setting $A:=h(L+R)$ gives
\begin{equation}
\label{eq:metdp_theta_form}
    I = h\int_0^1 e^{(1-\theta)A}\,e^{h(1-\theta)\ad_R}\bigl(N(t_n+h\theta)\bigr)\,d\theta.
\end{equation}

Expanding the exponential of $\ad_R$ as a power series,
\begin{equation*}
e^{h(1-\theta)\ad_R}=\sum_{j=0}^{\infty}\frac{h^j(1-\theta)^j}{j!}\,\ad_R^{\,j},
\end{equation*}
and inserting \eqref{eq:metdp_interp} into \eqref{eq:metdp_theta_form} yields the expansion
\begin{equation}
\label{eq:metdp_formal_series}
    I
    = h\sum_{m=0}^{p-1}\sum_{j=0}^{\infty} h^j\,C_{m,j}(A)\,
    \ad_R^{\,j}\!\left(\nabla^m N_n\right)
    \;+\;\mathcal O(h^{p+1}),
\end{equation}
where the coefficient operators admit the finite $\varphi$-representation
\begin{equation}
\label{eq:metdp_Cmj}
    C_{m,j}(A) = \frac{(-1)^m}{j!}\sum_{q=0}^{m+j} \alpha_{m,j,q}\,q!\,\varphi_{q+1}(A).
\end{equation}
and the coefficients $\alpha_{m,j,q}$ are defined as the coefficients of $\theta^q$ in the expansion of $(1-\theta)^j\binom{-\theta}{m}$ (see Appendix~\ref{app:METDp_deriv}).

\subsection{The METD\texorpdfstring{$p$}{p} update}
Along smooth solutions, one has $\nabla^m N_n=\mathcal O(h^m)$, so the term indexed by $(m,j)$ in \eqref{eq:metdp_formal_series} is of size $\mathcal O(h^{1+m+j})$ (precise bound given by Lemma~\ref{lem:diff-deriv} and the resulting local truncation error estimates are established in Section~\ref{sec:error}).
To obtain an order-$p$ method, we retain all contributions satisfying $m+j\le p-1$ and truncate the remainder.
Define
\begin{equation}
\label{eq:metdp_Ip}
    I_p := h\sum_{m=0}^{p-1}\sum_{j=0}^{p-1-m} h^j\,C_{m,j}(A)\,
    \ad_R^{\,j}\!\left(\nabla^m N_n\right),
    \qquad A=h(L+R).
\end{equation}
The resulting order-$p$ METD multistep scheme is
\begin{equation}
\label{eq:metdp_scheme}
    Q_{n+1}=e^{hL}Q_ne^{hR}+I_p.
\end{equation}
For $p=1$ and $p=2$, \eqref{eq:metdp_scheme} reduces to METD1 \eqref{eq:etd1-eqn} and METD2 \eqref{eq:etd2-eqn}, respectively, recovering the previously derived schemes. The complete integration procedure is summarized in Algorithm~\ref{alg:metdp} and the construction of the coefficient matrices $C_{m,j}$ is detailed in Algorithm~\ref{alg:build_cmj_algorithmic}.

\begin{algorithm}[t]
\caption{Arbitrary order METD$p$}
\label{alg:metdp}
\begin{algorithmic}[1]
\REQUIRE $Q_0$, step size $h$, order $p$, number of steps $N_t$, matrices $L,R$, nonlinearity $N(\cdot)$
\STATE Precompute $E_L=e^{hL}$ and $E_R=e^{hR}$
\STATE Set $A=h(L+R)$ and compute $\varphi_k(A)$ for $k=1,\ldots,p$
\STATE Build $C_{m,j}(A)$ for $0\le m\le p-1$, $0\le j\le p-1-m$
\FOR{$n=0,\ldots,N_t-1$}
    \STATE $N_n \leftarrow N(Q_n)$
    \STATE Build $\nabla^0 N_n,\nabla^1 N_n,\ldots,\nabla^{p-1}N_n$ from history
    \STATE $Q_{n+1} \leftarrow E_L Q_n E_R$
    \FOR{$m=0,\ldots,p-1$}
        \STATE $Y \leftarrow \nabla^m N_n$
        \FOR{$j=0,\ldots,p-1-m$}
            \STATE $Q_{n+1} \leftarrow Q_{n+1} + h^{j+1} C_{m,j}(A)\,Y$
            \STATE $Y \leftarrow YR-RY$
        \ENDFOR
    \ENDFOR
\ENDFOR
\RETURN $Q_{N_t}$
\end{algorithmic}
\end{algorithm}

\begin{rmk}
\label{rmk:metd2_consistency}
Specializing the METD$p$ construction to $p=2$ produces a single commutator correction term involving
$\varphi_k\,\big(h(L+R)\big)$. In that commutator contribution one may replace $e^{s(L+R)}$ by $e^{sL}$
without changing the method up to $\mathcal{O}(h^3)$. Hence replacing $\varphi_k\,\big(h(L+R)\big)$ by $\varphi_k(hL)$ in the
commutator term changes the one-step update only by $\mathcal{O}(h^3)$, so second-order accuracy is
preserved, recovering the METD2 form in Section~\ref{sec:metd2}. For $p\ge 3$ we instead use the exact
conjugation $e^{-sR}(\cdot)e^{sR}=e^{s\ad_R}(\cdot)$ and work with nested commutators $\ad_R^{\,j}$.
\end{rmk}

\section{Error analysis}
\label{sec:error}

Before the subsequent analysis, we adopt the following definition

\begin{defi}\label{def:1}
    Let $\mathcal{B}$ be a Banach algebra, i.e., a triple $(\mathcal{B}, \ast, ||\cdot||)$ where $(\mathcal{B}, \ast)$ is a unital associative algebra, and $(\mathcal{B},||\cdot||)$ is a (real or complex) Banach space, where the norm $||\cdot||$ is compatible with the multiplication $||x \ast y|| \leq ||x||\,||y||$ for any $x,y \in \mathcal{B}$. The commutator is then defined as $[x,y] = x \ast y - y \ast x$.
\end{defi}

In practice, we discretize our domain such that we will have $\mathcal{B}$ being the space of real (or complex) $n \times n$ matrices. We also make the following assumptions

\begin{asp}
    \label{asm:ass1}
    We assume that $L,R\in \mathcal B$ are bounded linear operators. Further, we assume that $L$ is an infinitesimal generator of an analytic semigroup $e^{tL}$ and the same for $R$ and $L+R$.
\end{asp}
This assumption implies that there exist constants $C>0$ and $\omega \in \mathbb{R}$ such that
\begin{equation}
    ||e^{tL}|| \leq C e^{\omega t}\, ,\qquad t \geq 0. 
\end{equation}
In particular, we have that $\varphi_k$ in~(\ref{eq:varphi_def}) for $k\geq1$ are bounded operators. 

\begin{asp}
    \label{asm:ass2}
    We assume that the nonlinearity $N : \mathcal{B} \times [0,T] \to \mathcal{B}$ is globally Lipschitz continuous on a strip along the exact solution. Moreover, along the exact solution the map $t\mapsto N(Q(t),t)$ is $C^{p}$ and its derivatives up to order $p$ are uniformly bounded on $[0,T]$.
\end{asp}

\subsection{Local defect of METD\texorpdfstring{$p$}{p}}

Let $Q(t)$ denote the exact solution of \eqref{eq:ODE} and set $t_n=nh$, $h=\Delta t$.
Recall the METD$p$ update \eqref{eq:metdp_scheme}, and define the one-step defect $\delta_{n+1}$ by inserting the exact solution into the numerical update
\begin{equation*}
Q(t_{n+1}) = e^{hL}Q(t_n)e^{hR} + I_p \;+\; \delta_{n+1},
\end{equation*}
where $I_p$ is the truncated approximation of the nonlinear integral \eqref{eq:metdp_Ip}.

\begin{thm}
\label{thm:metdp_local}
Assume \ref{asm:ass1}--\ref{asm:ass2} and $[L,R]=0$. Then there exists a constant $C$ (independent of $h$ and $n$ with $t_{n+1}\le T$) such that the defect satisfies
\begin{equation*}
\|\delta_{n+1}\|\le C\,h^{p+1}.
\end{equation*}
In particular, METD$p$ is consistent of order $p$.
\end{thm}

\begin{pf}
Write the exact nonlinear contribution over one step as
\begin{equation*}
I := \int_0^h e^{(h-\tau)L}\,N(Q(t_n+\tau),t_n+\tau)\,e^{(h-\tau)R}\,d\tau,
\end{equation*}
so that $\delta_{n+1}= I - I_p$.
We split the defect into an interpolation part and a commutator-truncation part,
\begin{equation*}
\delta_{n+1} = E_{\mathrm{interp}} + E_{\mathrm{comm}}.
\end{equation*}

\paragraph{Interpolation error ($E_{\mathrm{interp}}$).}
Define the Banach-space valued function
\begin{equation}
G(\tau):=N(Q(t_n+\tau),t_n+\tau),\qquad 0\le \tau\le h.
\end{equation}
Let $P_{p-1}(\tau)$ be the degree $(p-1)$ Newton backward interpolant of $G$ based on the grid values
$G(0)=N_n,\,G(-h)=N_{n-1},\dots,G(-(p-1)h)=N_{n-p+1}$.
Then
\begin{equation}
G(\tau)=P_{p-1}(\tau)+R_p(\tau),
\end{equation}
where the standard interpolation remainder $R_p(\tau)$ for equally spaced nodes yields the bound
\begin{equation}
\|R_p(\tau)\|\le C\,h^p\;\sup_{\xi\in[t_{n-p+1},\,t_{n+1}]}\left\|\frac{d^p}{dt^p}N(Q(t),t)\Big|_{t=\xi}\right\|
\le C\,h^p,
\qquad 0\le \tau\le h,
\end{equation}
using Assumption~\ref{asm:ass2}.
Inserting $R_p$ into the variation-of-constants integral gives
\begin{equation}
E_{\mathrm{interp}}
= \int_0^h e^{(h-\tau)L}\,R_p(\tau)\,e^{(h-\tau)R}\,d\tau.
\end{equation}
Using the semigroup bounds from Assumption~\ref{asm:ass1},
\begin{align*}
\|E_{\mathrm{interp}}\|
&\le \int_0^h \|e^{(h-\tau)L}\|\,\|R_p(\tau)\|\,\|e^{(h-\tau)R}\|\,d\tau \\
&\le C \int_0^h e^{2\omega(h-\tau)}\, (C h^p)\,d\tau
\le C\,h^{p+1}.
\end{align*}

\paragraph{Commutator truncation error ($E_{\mathrm{comm}}$).}
In the METD$p$ construction, after inserting the backward-difference expansion
\eqref{eq:metdp_interp}, the operator $e^{(h-\tau)\ad_R}$ is expanded in a series and truncated.
For each $m\in\{0,\dots,p-1\}$, METD$p$ keeps commutator depths $j\le p-1-m$, so the discarded tail starts at $j=p-m$.
Thus the truncation remainder has the form
\begin{equation}
E_{\mathrm{comm}}
= \sum_{m=0}^{p-1}\int_0^h e^{(h-\tau)(L+R)}
\left[\sum_{j=p-m}^{\infty}\frac{(h-\tau)^j}{j!}\,\ad_R^{\,j}(\nabla^m N_n)\right]\,d\tau,
\end{equation}
up to uniformly bounded scalar factors from $(1-\theta)^j\binom{-\theta}{m}$, which we absorb into the constant $C$ of theorem~\ref{thm:metdp_local}.
We bound nested commutators in a Banach algebra by
\begin{equation*}
\|[X,R]\|\le \|XR\|+\|RX\|\le 2\|X\|\,\|R\|,
\qquad\Rightarrow\qquad
\|\ad_R^{\,j}(X)\|\le 2^j\|X\|\,\|R\|^j.
\end{equation*}
Moreover, along the exact solution, Lemma~\ref{lem:diff-deriv} implies the backward-difference scaling
$\|\nabla^m N_n\|\le C h^m$ for $m\le p-1$.
Hence, for $0\le \tau\le h$,
\begin{align*}
\left\|\sum_{j=p-m}^{\infty}\frac{(h-\tau)^j}{j!}\,\ad_R^{\,j}(\nabla^m N_n)\right\|
&\le \sum_{j=p-m}^{\infty}\frac{h^j}{j!}\,2^j\,\|\nabla^m N_n\|\,\|R\|^j \\
&\le C h^m\sum_{j=p-m}^{\infty}\frac{(2h\|R\|)^j}{j!}
\;\le\; C h^m\,(2h\|R\|)^{p-m},
\end{align*}
where in the last step we used a standard bound on the tail of the exponential series (absorbing constants into $C$).
Therefore, the integrand is $\mathcal{O}(h^p)$, and using again Assumption~\ref{asm:ass1},
\begin{equation*}
\|E_{\mathrm{comm}}\|
\le \int_0^h \|e^{(h-\tau)(L+R)}\|\, (C h^p)\,d\tau
\le C\int_0^h e^{\omega(h-\tau)}\,h^p\,d\tau
\le C\,h^{p+1}.
\end{equation*}

We conclude
\begin{equation*}
    \|\delta_{n+1}\|
\le \|E_{\mathrm{interp}}\|+\|E_{\mathrm{comm}}\|
\le C h^{p+1},
\end{equation*}

which proves the claim.
\end{pf}

\subsection{Global convergence of METD\texorpdfstring{$p$}{p}}

\begin{thm}
\label{thm:metdp_global}
Assume \ref{asm:ass1}--\ref{asm:ass2} and $[L,R]=0$.
Let $Q(t)$ denote the exact solution of \eqref{eq:ODE} on $[0,T]$.
Let $\{Q_n\}_{n=0}^{M}$ be generated by the METD$p$ scheme \eqref{eq:metdp_scheme} with step size $h=T/M$.

Assume the starting values satisfy
\begin{equation}\label{eq:startup_assumption}
\max_{0\le n\le p-1}\|Q(t_n)-Q_n\|\le C_0\,h^p
\end{equation}
for some constant $C_0\ge 0$ (e.g.\ exact starting values, or produced by a one-step method of order $p$).

Then there exists a constant $C>0$, independent of $h$ (for $h\le h_0$) but depending on $T$, semigroup bounds for $e^{tL}$ and $e^{tR}$, the Lipschitz constant of $N$, and bounds on $t\mapsto N(Q(t),t)$ and its derivatives up to order $p$, such that
\begin{equation}\label{eq:global_error}
\max_{0\le n\le M}\|Q(t_n)-Q_n\|\le C\,h^{p}.
\end{equation}
\end{thm}

\begin{pf}
Let $\Phi_h:\mathcal{B}^p\to\mathcal{B}$ denote the METD$p$ update map, so that
\begin{equation*}
Q_{n+1}=\Phi_h(Q_n,\ldots,Q_{n-p+1})
      = e^{hL}Q_ne^{hR} + I_p(Q_n,\ldots,Q_{n-p+1}),
\end{equation*}
where $I_p$ is defined in \eqref{eq:metdp_Ip}. Define the global error
\begin{equation*}
e_n := Q(t_n)-Q_n,
\end{equation*}
and define the local defect by inserting the exact values into the update,
\begin{equation*}
\delta_{n+1}
:= Q(t_{n+1})-\Phi_h\bigl(Q(t_n),\ldots,Q(t_{n-p+1})\bigr).
\end{equation*}
By Theorem~\ref{thm:metdp_local}, there exists $C_{\mathrm{local}}>0$ such that
\begin{equation}\label{eq:glob_loc_def}
\|\delta_{n+1}\|\le C_{\mathrm{local}}\,h^{p+1}
\qquad\text{for all }n\text{ with }t_{n+1}\le T.
\end{equation}

Subtracting the numerical update from the exact update and using the definition of $\Phi_h$ we obtain the error recursion
\begin{equation}\label{eq:err_decomp}
e_{n+1}
= e^{hL}e_ne^{hR}
+ \Bigl(I_p(Q(t_n),\ldots,Q(t_{n-p+1}))-I_p(Q_n,\ldots,Q_{n-p+1})\Bigr)
+ \delta_{n+1}.
\end{equation}

Next, we claim that there exists a constant $K>0$, independent of $h$ and $n$ (for $t_{n+1}\le T$),
such that
\begin{equation}\label{eq:Ip_Lip}
\Bigl\|I_p(Q(t_n),\ldots,Q(t_{n-p+1}))-I_p(Q_n,\ldots,Q_{n-p+1})\Bigr\|
\le hK\,\max_{0\le \ell\le p-1}\|e_{n-\ell}\|.
\end{equation}
To see this, note from \eqref{eq:metdp_Ip} that $I_p$ is a finite sum (for fixed $p$) of terms of the form
\begin{equation*}
h^{1+j}\,C_{m,j}(A)\,\ad_R^{\,j}\bigl(\nabla^m N_n\bigr),
\qquad 0\le m\le p-1,\quad 0\le j\le p-1-m.
\end{equation*}
For $h\in(0,T]$, the operator coefficients $C_{m,j}(A)$ are bounded because they are finite linear
combinations of $\varphi$-functions evaluated at $A=h(L+R)$, and $\varphi_k$ are bounded operators under
Assumption~\ref{asm:ass1}. Moreover, $\ad_R$ is bounded since $R$ is bounded, hence $\ad_R^{\,j}$ is bounded for each
fixed $j$. Finally, by Assumption~\ref{asm:ass2}, $N(\cdot,t)$ is globally Lipschitz with constant $L_N$, so for each $\ell$
\begin{equation*}
\|N(Q(t_{n-\ell}),t_{n-\ell})-N(Q_{n-\ell},t_{n-\ell})\|\le L_N\|e_{n-\ell}\|.
\end{equation*}
Since $\nabla^m$ is a fixed finite linear combination of the past values $N_{n-\ell}$ (with coefficients depending only on $m$),
the difference $\nabla^m N(Q(t_\cdot))-\nabla^m N(Q_\cdot)$ is bounded by a constant (depending only on $p$) times
$\max_{0\le\ell\le p-1}\|e_{n-\ell}\|$. Collecting the finitely many bounded operator norms and combinatorial constants into $K$
yields \eqref{eq:Ip_Lip}.

Now, let us define 
\begin{equation*}
F_{n+1}
:= \Bigl(I_p(Q(t_n),\ldots,Q(t_{n-p+1}))-I_p(Q_n,\ldots,Q_{n-p+1})\Bigr)+\delta_{n+1}.
\end{equation*}
Then \eqref{eq:err_decomp} becomes the linear inhomogeneous recursion
\begin{equation}\label{eq:lin_inhom_rec}
e_{n+1}=e^{hL}e_ne^{hR}+F_{n+1}.
\end{equation}
Iterating \eqref{eq:lin_inhom_rec} gives the discrete variation-of-constants
\begin{equation*}
e_n
= e^{t_nL}e_0e^{t_nR}
 + \sum_{k=0}^{n-1} e^{(t_n-t_{k+1})L}\,F_{k+1}\,e^{(t_n-t_{k+1})R}.
\end{equation*}
Assuming exact initial data (or incorporating $e_0$ into the start-up assumption), we have $e_0=0$.
Using Assumption~\ref{asm:ass1}, for $0\le t_n-t_{k+1}\le T$,
\begin{equation*}
\|e^{(t_n-t_{k+1})L}\|\le C_L e^{\omega_L (t_n-t_{k+1})},\qquad
\|e^{(t_n-t_{k+1})R}\|\le C_R e^{\omega_R (t_n-t_{k+1})}.
\end{equation*}
Hence,
\begin{equation}\label{eq:en_sum_bound}
\|e_n\|
\le C_T\sum_{k=0}^{n-1}\|F_{k+1}\|,
\qquad
C_T:=C_LC_R\,e^{(\omega_L+\omega_R)T}.
\end{equation}
We define the maximum error
\begin{equation*}
   E_n:=\max_{0\le j\le n}\|e_j\|. 
\end{equation*}
Combining \eqref{eq:Ip_Lip} and \eqref{eq:glob_loc_def} yields, for all $k$ with $t_{k+1}\le T$,
\begin{equation*}
\|F_{k+1}\|\le hK\max_{0\le\ell\le p-1}\|e_{k-\ell}\| + C_{\mathrm{local}}h^{p+1}
\le hK\,E_k + C_{\mathrm{local}}h^{p+1}.
\end{equation*}
Inserting this into \eqref{eq:en_sum_bound} and using $nh\le T$ gives
\begin{align*}
\|e_n\|
&\le C_T\sum_{k=0}^{n-1}\bigl(hK E_k + C_{\mathrm{local}}h^{p+1}\bigr)
\le C_TK h\sum_{k=0}^{n-1}E_k + C_T C_{\mathrm{local}}\,T\,h^p.
\end{align*}
Taking the maximum over $0\le j\le n$ (so that $\|e_n\|\le E_n$) yields
\begin{equation}\label{eq:En_discrete}
E_n \le a h\sum_{k=0}^{n-1}E_k + b h^p,
\qquad
a:=C_TK,\quad b:=C_T C_{\mathrm{local}}\,T.
\end{equation}
For $n\ge p-1$, we can split the sum in \eqref{eq:En_discrete} as
\begin{equation*}
\sum_{k=0}^{n-1}E_k
\le \sum_{k=0}^{p-2}E_k + \sum_{k=p-1}^{n-1}E_k
\le (p-1)\max_{0\le j\le p-1}E_j + \sum_{k=p-1}^{n-1}E_k.
\end{equation*}
By the start-up assumption \eqref{eq:startup_assumption},
$\max_{0\le j\le p-1}E_j\le C_0 h^p$, so
\begin{equation*}
\sum_{k=0}^{n-1}E_k \le (p-1)C_0 h^p + \sum_{k=p-1}^{n-1}E_k.
\end{equation*}
Define $\widetilde E_n:=\max_{p-1\le j\le n}E_j$. Then from \eqref{eq:En_discrete} and the bound above,
for $n\ge p-1$,
\begin{equation*}
E_n \le a h\sum_{k=p-1}^{n-1}E_k + \bigl(a h (p-1)C_0 + b\bigr)h^p
\le a h\sum_{k=p-1}^{n-1}\widetilde E_k + \widetilde b\,h^p,
\end{equation*}
where $\widetilde b:=\bigl(a (p-1)C_0 + b\bigr)$.
Taking the maximum over $p-1\le j\le n$ gives
\begin{equation*}
\widetilde E_n \le a h\sum_{k=p-1}^{n-1}\widetilde E_k + \widetilde b\,h^p.
\end{equation*}
By the discrete Gr\"onwall inequality, this implies
\begin{equation*}
\widetilde E_n \le \widetilde b\,h^p e^{a(t_n-t_{p-1})}
\le \widetilde b\,e^{aT}\,h^p.
\end{equation*}
Since $E_n=\max\{\max_{0\le j\le p-1}E_j,\widetilde E_n\}$ and $\max_{0\le j\le p-1}E_j\le C_0 h^p$,
we conclude that for all $0\le n\le M$,
\begin{equation*}
E_n \le C\,h^p,
\qquad
C:=\max\{C_0,\widetilde b\,e^{aT}\}.
\end{equation*}
This proves \eqref{eq:global_error}.

\end{pf}

\section{Extension to Non-Commuting Operators}
\label{sec:extensions_non_comm}
We now consider the case $[L,R]\neq 0$. In the commuting derivation we used the identity
$e^{sL}e^{sR}=e^{s(L+R)}$, which no longer holds when $L$ and $R$ do not commute. A natural
replacement is obtained through the
Baker--Campbell--Hausdorff (BCH) Lie series.\\

\subsection{The Baker-Campbell-Hausdorff Series}

Let $\mathcal{B}$ be a (complex) Banach algebra with identity and let $L,R\in\mathcal{B}$.
For a step size $h>0$ define
\begin{equation*}
W_h := e^{hL}e^{hR}.
\end{equation*}
Whenever the principal logarithm of $W_h$ exists (equivalently, its spectrum does not intersect $(-\infty,0]$), we define
\begin{equation}
\label{eq:Zh_def}
Z_h := \log\,(W_h)=\log\,\bigl(e^{hL}e^{hR}\bigr),
\end{equation}
so that $e^{Z_h}=e^{hL}e^{hR}$. In the commuting case one has $Z_h=h(L+R)$.

For sufficiently small $h$, $Z_h$ admits the BCH expansion
\begin{equation}
\label{eq:BCH_series}
Z_h
= h(L+R)
+ \frac{h^2}{2}[L,R]
+ \frac{h^3}{12}\Bigl([L,[L,R]] + [R,[R,L]]\Bigr)
+ \mathcal{O}(h^4)\,.
\end{equation}

Moreover, the BCH series converges under a smallness condition

\begin{thm}[Rossmann \cite{rossmann2006lie}]
\label{thm:BCH_convergence}
If $X,Y\in\mathcal{B}$ satisfy $\|X\|+\|Y\|<\log(2)/2$, then the BCH series converges
absolutely and its sum equals $\log(e^Xe^Y)$.
\end{thm}

In our setting $X=hL$ and $Y=hR$, so Theorem~\ref{thm:BCH_convergence} guarantees convergence
whenever $h(\|L\|+\|R\|)<\log(2)/2$. Outside this regime, one may still compute $Z_h$
directly using a matrix-logarithm algorithm whenever a suitable logarithm exists. \\

\subsection{METD for $[L,R] \neq 0$}
Starting from the variation-of-constants formula \eqref{eq:etd-int}, we use the conjugation
identity $e^{s\,\ad_R}(X)=e^{-sR}Xe^{sR}$ to rewrite
\begin{equation}
\label{eq:conjugation_BCH}
e^{sL} N e^{sR}
= \bigl(e^{sL}e^{sR}\bigr)\,\bigl(e^{-sR}Ne^{sR}\bigr)
= e^{Z_s}\,e^{s\,\ad_R}(N),
\end{equation}
where $Z_s:=\log(e^{sL}e^{sR})$ whenever the logarithm is defined.
With the change of variables $s=h-\tau$, the nonlinear contribution over one step becomes
\begin{equation}
\label{eq:exact_integral_BCH}
I = \int_0^h e^{Z_s}\,e^{s\,\ad_R}\,\bigl(N(t_{n+1}-s)\bigr)\,ds.
\end{equation}

Unlike the commuting case, $Z_s$ is not linear in $s$ in general. To obtain a practical
fixed-step method we approximate $Z_s$ by a linear interpolation between the endpoints,
\begin{equation}
\label{eq:linear_freezing}
Z_s \approx \frac{s}{h}Z_h\,.
\end{equation}
Under the validity of the BCH expansion for small step sizes, this replacement satisfies
$\|Z_s-\tfrac{s}{h}Z_h\|=\mathcal O(h^2)$ uniformly for $s\in[0,h]$, and thus it is
compatible with at most second-order accuracy unless higher-order approximations to $Z_s$
are employed.

When $Z_h$ is not computed directly via a matrix-logarithm routine, we approximate it
using a truncated BCH expansion,
\begin{equation}
\label{eq:Zh_trunc}
Z_h \approx Z_h^{(n)}:=h(L+R)+\frac{h^2}{2}[L,R]+\cdots,\qquad (\text{commutators up to depth }n),
\end{equation}
where $n\ge 1$ denotes the truncation depth in nested commutators. We denote the resulting
schemes by METD$j$\_BCH$n$, where $j$ is the temporal order of the METD update and $n$ is the
BCH truncation depth.

For example, METD1-BCH can be obtained by using \eqref{eq:linear_freezing} which yields
\begin{equation}
\label{eq:bch_metd1}
Q_{n+1} = e^{hL}Q_ne^{hR}
+ h\,\varphi_1\bigl(Z_h\bigr)\,N_n
\;\;\; +\;\; \mathcal O(h^2),
\end{equation}
where $\varphi_1$ is defined as in \eqref{eq:varphi_def} and is evaluated at $Z_h$ (or at
$Z_h^{(n)}$ when using BCH truncation).

\begin{rmk}
\label{rmk:order_limit}
The linear freezing approximation \eqref{eq:linear_freezing} introduces an $\mathcal{O}(h^2)$ 
error in the integrand, which becomes $\mathcal{O}(h^3)$ after integration. Consequently, 
METD$j$-BCH$n$ schemes are limited to at most second-order accuracy regardless of $j$ or $n$. 
Achieving higher-order convergence in the non-commuting case would require retaining the 
explicit $s$-dependence in the BCH expansion of $Z_s$ directly.
\end{rmk}

\section{Numerical Examples}
In this section, we compute and test the METD$p$ schemes with differential Lyapunov and Sylvester equations, and with a stiff nonlinear PDE-derived matrix system (Allen–Cahn). Our main goal is to confirm the global truncation error of the schemes and show that they are applicable to real-world, high-dimensional complex systems from PDEs to machine learning. We start with a simple test and gradually increase complexity and generality.

\subsection{A low-dimensional Lyapunov equation}
\label{sec:num_lyap}

\begin{figure}
    \centering
    \includegraphics[width=0.6\linewidth]{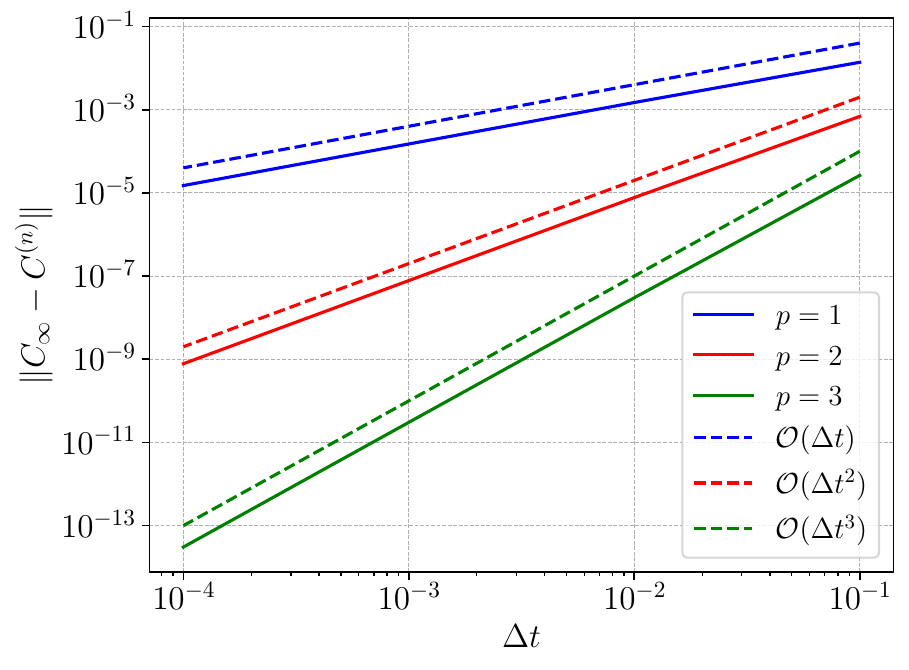}
    \caption{Error at $T=10$ for the Lyapunov equation~\eqref{eq:lyap} measured as $\|C_\infty - C^{(n)}\|$ versus $\Delta t$, where $C_\infty$ solves the stationary Lyapunov equation. The observed slopes confirm the expected global orders for METD$p$ (shown here for $p=1,2,3$).}
    \label{fig:etd1_toy2}
\end{figure}

As the simplest non-trivial example of a matrix evolution equation, we consider the dynamical Lyapunov system
\begin{equation}
    \label{eq:lyap}
    \dot C(t) = A\, C(t) + C(t) \, A^T + \Sigma,
\end{equation}
where $C(t)\in\mathbb{R}^{n \times n}$, $\Sigma \in \mathbb{R}^{n \times n}$ symmetric positive definite, and $A \in \mathbb{R}^{n \times n}$ normal ($AA^T=A^TA$) and negative definite. While this example is merely affine-linear, with $N=\Sigma$ independent of $C(t)$, it is still non-trivial for METD$p$ because commutator
corrections with $R=A^T$ need not vanish (in general $[\Sigma,A^T]\neq 0$). Moreover, because $N$ is constant in time, all backward differences
$\nabla^m N_n$ vanish for $m\ge 1$. As a result, this test cleanly isolates the effect of truncating the commutator contributions at increasing depth
as $p$ grows. 

We solve the system~(\ref{eq:lyap}) until time $T=10$, where it is close enough to the solution $C_\infty$ of the stationary Lyapunov equation that we can take $\|C(T)-C_\infty\|$ as the error (i.e. error is dominated by time discretization error). Solving this system for $n=2$, we see in figure \ref{fig:etd1_toy2} that the observed slopes match the expected global orders: METD1 is first-order, METD2 is second-order, and METD3 is third order. Computationally, for dense matrices, the dominant cost is forming matrix exponentials (and $\varphi$-functions), giving $\mathcal{O}(n^3)$ complexity per step.

We next consider a system where the nonlinear term is truly nonlinear.

\subsection{Allen--Cahn equation}
\label{sec:allen_cahn} 

Let us consider the Allen--Cahn equation on a torus $(x,y)\in\TT^2_{2\pi} = [0,2\pi]^2$,
\begin{equation}
  \label{eq:allen_cahn_pde}
  \partial_t f = \varepsilon\,\Delta f + f - f^3,
\end{equation}
with $\varepsilon = 0.1$ and periodic boundary conditions. This problem, which
models phase separation ~\cite{allen1972ground},
has also been used as a benchmark in~\cite{rodgers2023implicit,carrel2023projected}.
Similarly to both papers cited above, we have the initial profile
\begin{equation}
  \label{eq:ac_ic}
  f_0(x,y) =
  \frac{\bigl(e^{-\tan^2 x} + e^{-\tan^2 y}\bigr)\,\sin x\,\sin y}
       {1 + \exp\,\bigl(|\csc(-x/2)|\bigr) + \exp\,\bigl(|\csc(-y/2)|\bigr)} \,.
\end{equation}

A fourth-order finite-difference spatial discretization on a $256\times 256$ grid yields
the matrix ODE
\begin{equation}
  \label{eq:allen_cahn_matrix}
  \dot X = AX + XA + X - X^{\circ 3},
  \qquad X(0) = X_0,
\end{equation}
where $A\in\RR^{256\times 256}$ is the one-dimensional discrete Laplacian scaled by
$\varepsilon$, and $X^{\circ 3}$ denotes the element-wise (Hadamard) cube, i.e.\
$(X^{\circ 3})_{ij} = X_{ij}^3$. In the notation of~\eqref{eq:ODE}, we have
$L = R = A$ and $N(X) = X - X^{\circ 3}$. Because the
Laplacian eigenvalues scale as $\varepsilon k^2$, the linear part of \eqref{eq:allen_cahn_matrix} is severely stiff, limiting standard explicit methods to very small step sizes.

The reference trajectory is computed with SciPy's vectorized variable step RK45 solver
(with absolute tolerance set to $10^{-14}$). This reference provides
both ground truth for computing global errors and the start-up values
$Q_0,\dots,Q_{p-2}$ needed for the METD$p$ multistep schemes.

\begin{figure}
    \centering
    \includegraphics[width=\linewidth]{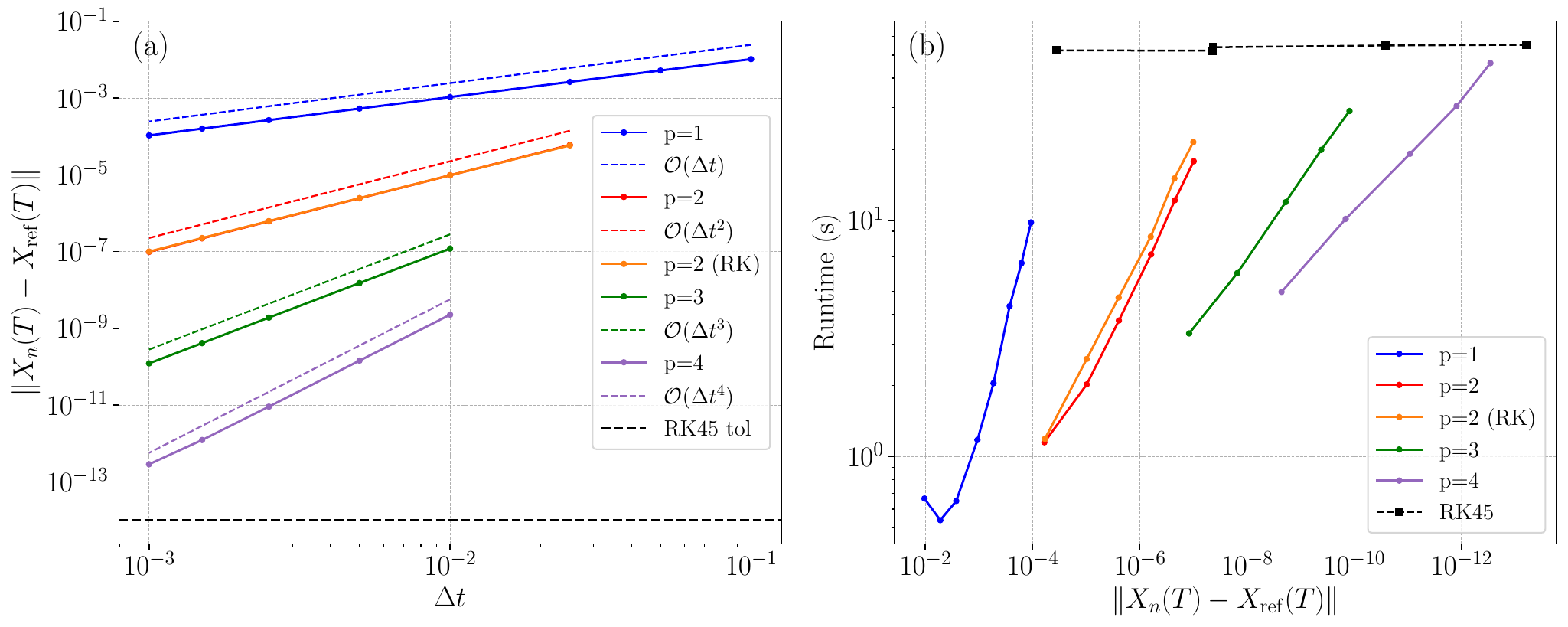}
    \caption{Allen--Cahn equation~\eqref{eq:allen_cahn_matrix} with
    $\varepsilon=0.1$ on a $256\times 256$ grid.
    (a)~The observed slopes confirm the expected convergence orders for
    METD$p$ with $p=1,2,3,4$ and the Runge--Kutta variant METD2RK.
    For $\Delta t < 10^{-3}$ the METD4 curve levels off at the numerical
    error floor.
    (b)~Each METD marker shows the error achieved across each step
    size; the corresponding RK45 marker shows the runtime required to
    reach the same accuracy. RK45 runtime is roughly constant because
    its step size is stability-limited on this stiff system.}
    \label{fig:allen_cahn}
\end{figure}

We benchmark METD1, METD2, METD2RK, METD3 and METD4 over a set of time steps $\Delta t\in\{10^{-1},\dots,10^{-3}\}$ up to $T=14$.
The left panel of figure~\ref{fig:allen_cahn} reports the global error
$\|X_n(T)-X_{\mathrm{ref}}(T)\|$ as a function of the step size $\Delta t$.
Every METD variant exhibits its theoretical convergence rate on this stiff
PDE. The higher-order methods (METD2 and above) become unstable at the largest step
sizes ($\Delta t \geq 0.05$), reflecting the practical stability boundary of
the explicit treatment of the nonlinearity. METD1, by contrast, remains stable
even at $\Delta t = 0.1$ and exhibits clean first-order convergence. 

To compare runtimes fairly, we match each METD method against RK45 at the
same accuracy. For each METD variant, we take the minimum relative
error $e_m^\star$ achieved across all tested step sizes, then run the
vectorized RK45 solver with tolerances set to $e_m^\star$ and record its runtime and achieved error. The right panel of figure~\ref{fig:allen_cahn} shows the result. The step size of RK45 on this system is dictated by stability rather than accuracy; its runtime is essentially constant ($\approx 30$\,s) regardless of the requested tolerance. METD, by contrast,
decouples stability from step size through the exact linear propagator
$e^{hA}(\cdot)e^{hA}$, so its cost scales directly with the number of steps
taken. At low-to-moderate accuracy this gives substantial speedups, METD1
is roughly $5\times$ faster than RK45 and METD2 roughly $2.5\times$
faster, while at the highest accuracies the runtimes become comparable. The key advantage is that by removing the stiff linear part from the
stability constraint, METD allows the step size to be dictated by the
nonlinearity rather than by the Laplacian eigenvalues, enabling a
meaningful accuracy-cost tradeoff over a wide range of step sizes.

\subsubsection{Comparison with alternative approaches}
\label{sec:comparison_alternatives}
\begin{table}[t]
\centering
\caption{Allen--Cahn comparison with step-truncation (ST) midpoint methods from \cite{rodgers2023implicit}. Runtimes for ST methods are taken from \cite{rodgers2023implicit} and may depend on implementation and hardware; comparisons should therefore be interpreted at the order-of-magnitude level.
}
\label{tab:allen_cahn_comparison}
\begin{tabular}{lcccc}
\hline
Method  & $\Delta t$ & $\|X_n(T)-X_{\mathrm{ref}}(T)\|$ & Runtime (s) \\
\hline
ST midpoint (explicit) \cite{rodgers2023implicit} & $10^{-3}$ & unstable & -- \\

ST midpoint (implicit) \cite{rodgers2023implicit} & $10^{-2}$ & $6.7\times 10^{-6}$ & 38.03 \\

ST midpoint (implicit) \cite{rodgers2023implicit} & $10^{-1}$ & $2.6\times 10^{-2}$ & 5.31 \\

METD1 & $10^{-1}$ & $1.0\times 10^{-2}$ & 0.67 \\

METD2 & $10^{-2}$ & $9.7\times 10^{-6}$ & 2.02 \\

METD4 & $10^{-2}$ & $2.3\times 10^{-9}$ & 4.98 \\
\hline
\end{tabular}
\end{table}
It is instructive to compare with the implicit step-truncation (ST) methods of
Rodgers and Venturi~\cite{rodgers2023implicit}, who solve the same
Allen--Cahn equation ($\varepsilon=0.1$, periodic domain, comparable grid). Their approach achieves unconditional stability by formulating each time step as a root-finding problem. METD achieves stability through a fundamentally different mechanism; the
stiff linear propagator is absorbed exactly into the matrix exponentials
$e^{hA}(\cdot)e^{hA}$, and the remaining nonlinearity is integrated
explicitly. This requires no iterative solves, no Krylov subspace
construction, and no rank management. After precomputing the exponential/$\varphi$ operators for a given $\Delta t$, each step simply consists of matrix multiplications and pointwise nonlinear evaluations. 
We compare METD with ST methods in table~\ref{tab:allen_cahn_comparison} which highlights two points. First, the explicit ST midpoint method becomes unstable at much smaller step sizes, whereas METD1 remains stable at $\Delta t=10^{-1}$. Second, at matched step size $\Delta t=10^{-2}$, METD2 achieves comparable accuracy to implicit ST midpoint at substantially lower runtime, while METD4 attains several additional digits of accuracy with only a modest increase in cost. At $\Delta t=10^{-1}$, METD1 is both significantly faster and also more accurate than implicit ST midpoint.

Carrel and Vandereycken~\cite{carrel2023projected} consider projected exponential
Runge--Kutta (PERK) methods within a dynamical low-rank approximation (DLRA)
framework on the Allen--Cahn system (in a less stiff parameter regime with $\varepsilon=0.01$). For Allen--Cahn, the elementwise cubic term is challenging because $\operatorname{rank}(X^{\circ 3}) \le \operatorname{rank}(X)^3$ so nonlinear evaluations can substantially increase intermediate rank. PERK/DLRA controls this growth through tangent-space projection and subsequent truncation to a prescribed rank $r$. In our numerical tests of PERK, increasing the target rank to a moderate size (e.g., to $r=20$) led to numerical difficulties (including SVD-truncation failures), whereas the dense METD discretization requires no rank truncation and is stable for large $\Delta t$. The full $n\times n$ state is propagated directly, and the nonlinearity is evaluated without projection, and no SVD-based compression is needed at any stage.

One may ask also whether the matrix structure is necessary at all. Vectorizing
\eqref{eq:allen_cahn_matrix} yields
\begin{equation*}
\dot{\hat Q}=\mathcal L \hat Q+\mathcal N(\hat Q),\qquad
\hat Q\in\mathbb R^{n^2},
\end{equation*}
with $\mathcal L = I\otimes A + A\otimes I \in \mathbb R^{n^2\times n^2}$. Our RK45 reference indeed integrates this vectorized system.
However, a direct scalar ETD discretization on the vectorized problem would
involve $\exp(h\mathcal L)$ or $\varphi_k(h\mathcal L)$ on an
$n^2\times n^2$ operator, which is $\mathcal O(n^6)$ in dense arithmetic. Matrix-free Krylov variants can avoid forming $\mathcal L$, but still work in
dimension $n^2$ and require repeated operator actions and basis storage. In contrast, METD keeps the $n\times n$ form throughout. For a fixed
$\Delta t$, the exponential/$\varphi$ operators are precomputed once; each
time step then uses matrix multiplications,
with $\mathcal O(n^3)$ arithmetic and $\mathcal O(n^2)$ state storage.

\subsection{Turbulent fluctuations of atmospheric jets}
\label{Jets}
\begin{figure*}
    \centering
    \includegraphics[width=\textwidth,height=7cm]{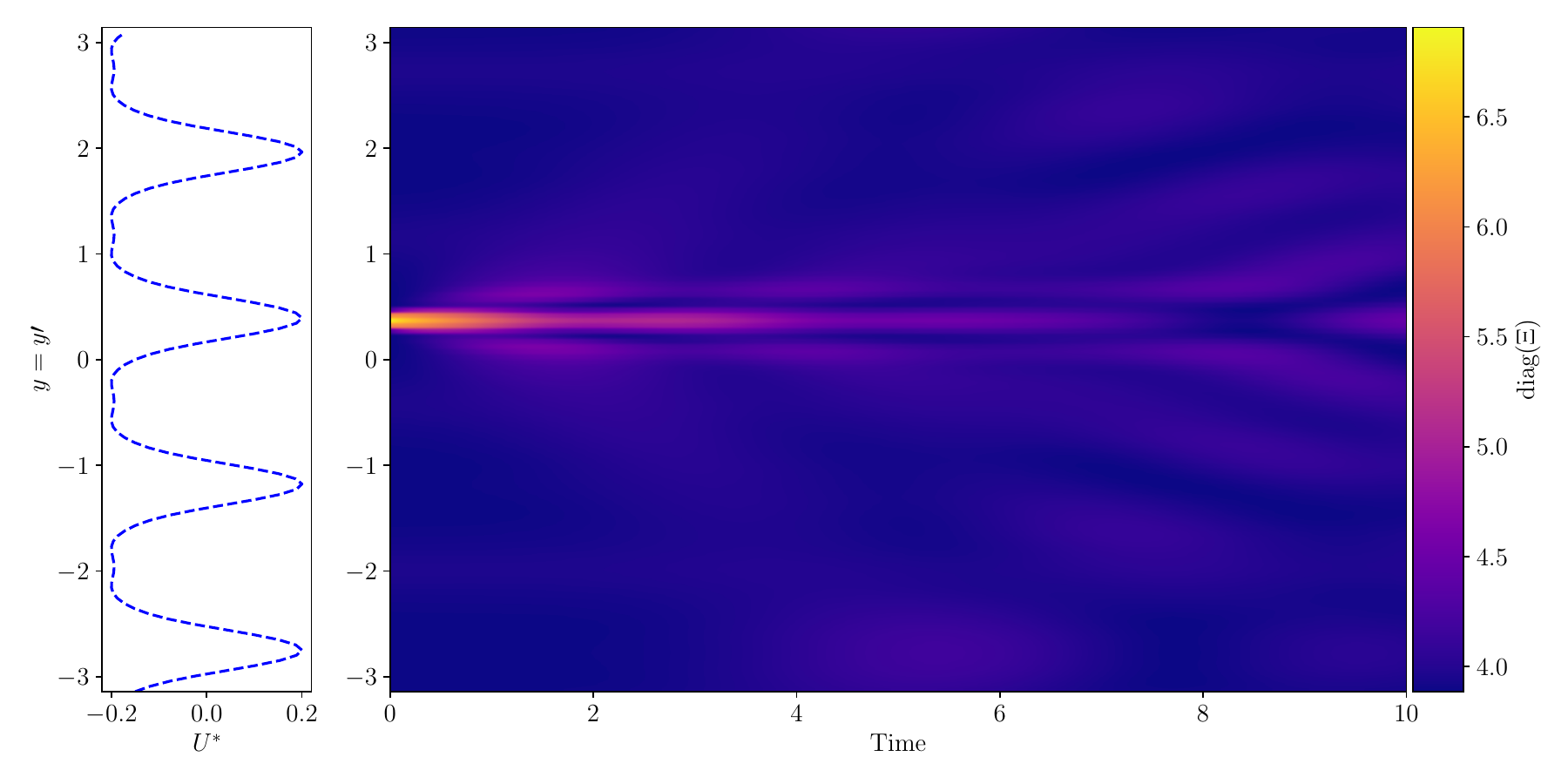}
    \caption{Perturbation on top of zonal jets decaying to a fixed point through the evolution of fluctuations via the differential Lyapunov equation \eqref{dynamicLyapjet}. The plot on the left shows the stable zonal jet profile $U^*$ with the stable $4$ jet configuration. Integrated using the METD2 scheme.}
    \label{fig:etd_pde}
\end{figure*}

To demonstrate performance at scale, we next consider a PDE-driven application from geophysical turbulence whose statistical closure yields a very large differential Lyapunov equation. We consider geostrophic turbulence in the atmosphere, for example in large gaseous planets, for simplicity modeled as a 2D reduction of the planetary atmosphere in a single layer. A simple model for atmospheric flow on a rotating planet is given by the 2-dimensional, stochastically forced, quasi-geostrophic equation in the $\beta$-plane \cite{quasi_beta}: For a periodic domain $(x,y)\in\TT^2_L = [0,L]^2$, the wind velocity vector field $\vec{v}=(u,v): \TT^2_L \mapsto \RR^2$, and its corresponding vorticity, $\omega = \nabla^\perp \vec{v}$, evolves according to
\begin{equation}
  \label{eq:QG}
  \partial_t \omega +\vec{v}\cdot\nabla \omega + \beta v = -\lambda
  \omega -\nu (-\Delta )^p\omega + \sqrt{\gamma}\eta\,.
\end{equation}
Here, the stream function $\psi:\TT^2_L\to\RR$ relates to the velocity via
$\vec v = e_z\times \nabla \psi$, and the vorticity fulfills $\omega = \Delta \psi$. Equation~(\ref{eq:QG})  includes a linear Ekman-damping term with friction coefficient $\lambda$,
and hyperviscosity of order $p$ with viscosity coefficient $\nu$. The
stochastic noise $\eta$ scales with amplitude $\gamma$, and has
covariance
\begin{equation}
  {\mathbb{E}}[\eta(x,y,t)\eta(x',y',t')] = \delta(t-t') C(x, y, x', y')\,,
\end{equation}
and the noise is homogeneous and isotropic, i.e.~$C(x, y, x', y')=\chi(|(x,y)-(x',y')|)$. Lastly, the $\beta$-term
models the Coriolis forces due to the planet's rotation, and determines the
$x$-direction to be the zonal direction. 

For a wide range of parameters, equation~(\ref{eq:QG}) exhibits the formation of large-scale coherent structures in the form of atmospheric jets, comparable to those observed on Jupiter or Saturn, and to a lesser degree Earth's jet stream. Defining $U(y)=\overline{u}(y) = \int u(x,y)\,dx$ the zonal velocity average, we can obtain an evolution equation for the zonal velocity, depending on higher-order turbulent vorticity fluctuations, by making the ansatz $u(x,y) \leftarrow U(y) + \sqrt{\alpha} u(x,y)$ and $v(x,y) \leftarrow \sqrt{\alpha} v(x,y)$~\cite{bouchet-nardini-tangarife:2013, bouchet-marston-tangarife:2018}. Considering $\Xi(x,y,x',y')={\mathbb{E}}[\omega(x,y)\omega(x',y')]$, the covariance of the two-point correlation function of vorticity fluctuations, the system we obtain is
\begin{equation}
    \begin{cases}
    \partial_t U &=-\alpha \overline{v\omega}   - \alpha U +\nu \partial_y^{2p} U\,\\
    \partial_t \Xi &= -(\Gamma(U) \, \Xi + \Xi \, \Gamma(U)^\dagger) + 2C\,,
    \end{cases}
\end{equation}
where $(\cdot) ^\dagger$ denotes the adjoint. In this approach, the nonlinear self-interaction of the fluctuation field is neglected, and we define
\begin{equation}
    \Gamma(U) := U\partial_x +(\partial_y^2 U-\beta)\partial_x \Delta^{-1} + \alpha + \nu(-\Delta)^p\,.
\end{equation}
Note that this system is commonly referred to as the second-order cumulant expansion (CE2) system \cite{ZonostrophicInstability, StructuralStabilityofTurbulentJets,PhysRevLett.110.224501}. It relates the slow evolution of atmospheric jets under the influence of fast turbulent fluctuations through the Reynolds stress~\cite{bouchet-simonnet:2009, bouchet-rolland-simonnet:2019}. We have a time scale separation in this system between the slow varying zonal jets and the fast turbulent fluctuations.  

By considering stable jet configurations $U^*(u)$, with $\partial_t U^* = 0$, perturbations of the turbulent fluctuations evolve via the Lyapunov system 
\begin{equation}
    \label{dynamicLyapjet}
    \partial_t \Xi = -(\Gamma(U^*) \, \Xi + \Xi \, \Gamma(U^*)^\dagger) + 2C\,.
\end{equation}
In this example, we discretize our problem into a grid of size $N_x \times N_y$, and use pseudo-spectral methods to solve our system. We choose $N_x = 64$ and $N_y=128$, resulting in a total of $n=N_x N_y = 8192$, so that $\Xi \in \RR^{8192\times 8192}$ to be solved by the METD scheme.

\textbf{Implementation details.} Our parameters are $\beta = 5$, $\alpha = 1 \times 10^{-3}$, $\nu = 1 \times 10^{-6} $ and $p=4$. We consider the $y-$direction in real space and take the $x-$direction in Fourier space. This allows us to decouple each Fourier mode and integrate equation~(\ref{dynamicLyapjet}) for each Fourier mode independently.

The particular difficulty of this problem is the fact that the $\Gamma$-operator includes a hyperviscosity term which is the Laplacian operator ($\Delta$) raised to the power $p$. This term dissipates fluctuations on the smallest scales (or highest modes) in our system which increases the range of resolved scales.  High orders of $p$ are often used in fluid-dynamics applications to effectively dampen fine scales of fluctuations without influencing larger scales, resulting in a broad band of physically meaningful length-scales not affected by viscosity. On the other hand, high orders of $p$ result in extremely harsh Courant–Friedrichs–Lewy (CFL) conditions, which for naive explicit solvers would require an extremely small time step. METD allows us to retain numerical stability even for very large time steps in this example (such as $\Delta t = 0.5$), due to the exact treatment of the linear $\Gamma$-operator. This improves integration performance of the system immensely.

As a concrete test case, we start the integration with an initial perturbation localized to a specific latitude, and can thus observe the propagation of the (variance of the) disturbance on the background of the atmospheric jet. The perturbation we choose is a simple Gaussian field centered at $\Xi_{y=y'}$ on a jet and we study the decay of this perturbation through the dynamic Lyapunov equation~(\ref{dynamicLyapjet}). Figure \ref{fig:etd_pde} shows the evolution of the perturbed diagonal of the correlation function. This example demonstrates the applicability of our matrix-ETD integration schemes in a complex setting.

\subsection{Extensions with Neural ODEs}

\subsubsection{Differential Sylvester equation}
We can generalize the METD scheme further to non-square problems. These systems are of the form shown in (\ref{eq:ODE}) except we now have $L \in  \mathbb{R}^{m \times m}$, $R \in  \mathbb{R}^{n \times n}$ with $Q,N \in  \mathbb{R}^{m \times n}$ where $m>n$, so the matrices $Q$ and $N(Q)$ are non-square. For the special case of $N$ constant, this system is known as the differential Sylvester equation and has a unique solution when the spectra of $L$ and $-R$ are disjoint \cite{behr_kyrlov}. The differential Lyapunov equation is a symmetric case ($L = R^T$) of the differential Sylvester equation. We cannot use the METD schemes presented in earlier sections without modification since we cannot consider the commutativity of matrices of different dimensions.

We can augment the operators such that we have a square system that considers all degrees of freedom. A valid augmentation is zero-padding the matrices to match the dimension of the largest matrix in the system. We would then have an augmented system
\begin{equation}
    \label{eq:ODE_augment}
    \dot{ \Tilde{Q} } = L\, \Tilde{Q} + \Tilde{Q} \, \Tilde{R} + \Tilde{N},
\end{equation}
where $\Tilde{Q}, \Tilde{N} \in \mathbb{R}^{m \times m}$ are zero-padded matrices with $m-n$ zero columns right-con\-cate\-nated to the original $Q$ and $N$ matrices and  $ \Tilde{R} = R \bigoplus \, 0_{(m-n)\times(m-n)}$ where $\bigoplus$ is the matrix direct sum. With this augmentation, we can consider the commutator $[ \Tilde{N}, \Tilde{R}]$ which allows us to approximate the integral in (\ref{eq:etd-int}). We only need to zero-pad when approximating the integral in (\ref{eq:etd-int}) and can easily remove the columns of zeros in our solution since the padded degrees of freedom do not ever interact with the original degrees of freedom. Of course, this strategy generalizes to an arbitrary non-square nonlinearity $N(Q)$.

Regarding the concrete implementation, we can use properties of the matrix exponential such that we do not ever compute the matrix exponential of the padded matrix $\Tilde{R}$ explicitly. The matrix exponential of a block zero matrix is the block identity matrix which means we can compute the matrix exponential of the smaller $R$ matrix and pad it with an identity matrix to avoid computing the matrix exponential of the larger $\Tilde{R}$ (shown by Lemma~\ref{matrix_exp_identity}). All matrix multiplications can be carried out with the original matrices since the padded degrees of freedom do not interact. These properties allow us to implement METD for the case $m>n$ with fewer operations than naively implementing METD for the augmented system (\ref{eq:ODE_augment}). 

\subsubsection{Continuous Graph Neural Networks}
The following example will make use of the generalizations discussed in the previous sections. Let us define a simple graph $G := (V,E)$ where $V$ is the set of vertices and $E \subseteq V \times V$ is the set of edges between vertices. Following Xhonneux et al.~\cite{pmlr-v119-xhonneux20a}, we can represent the graph structure using a (regularized) adjacency matrix 
\begin{equation}
    \Tilde{A}=: \frac{\alpha}{2} \left( \mathds{1} + D^{-\frac{1}{2}}\, A \, D^{-\frac{1}{2}} \right)\,,
\end{equation}

where $\alpha \in (0,1)$ is a hyperparameter,  $A \in \mathbb{R}^{ |V|\, \times \,|V|}$ is the typical adjacency matrix and $D_{ii} = \sum_{j}A_{ij}$ is the degree matrix. This modification, first introduced by Kipf et al.~\cite{kipf2016semi}, stabilizes graph learning by shifting the eigenvalues of the adjacency matrix to the interval $[0,\alpha]$.

In addition to the adjacency matrix, we have a finite set of features $F$, and a node feature matrix $X \in\mathbb{R}^{ |V|\, \times \, |F|} $ with $|F|$ being the number of node features. Our objective is to then learn a matrix of node representations $Q \in \mathbb{R}^{ |V| \, \times \, d} $ where $d$ is the dimension of the representation and the $i$-th row of $Q$ corresponds to the representation of the $i$-th node.

Graph neural networks (GNNs) \cite{kipf2016semi,velickovic2017graph,hamilton2017inductive} provide a robust framework for learning node representations by processing complex, non-Euclidean graph data. GNNs typically model the discrete dynamics of node representations using multiple propagation layers. In each layer, every node's representation is updated based on the representations of its neighboring nodes (message passing).

However, existing GNNs (e.g. GCN \cite{kipf2016semi}, GraphSage \cite{hamilton2017inductive}, GATs \cite{velickovic2017graph}) have been shown to suffer from a phenomenon called over-smoothing \cite{li2018deeper}. Over-smoothing occurs when node representations in a graph become increasingly similar through multiple layers of message passing, leading to a loss of discriminative power and ultimately, poor performance. An interesting model by Xhonneux et al. \cite{pmlr-v119-xhonneux20a} tries to mitigate the over-smoothing problem by taking inspiration from Neural ODEs \cite{chen2018neural} and adapting the discrete dynamics to a continuous space. We therefore have a continuous dynamical system to model the continuous dynamics on node representations. Moreover, it has been shown that this approach does not experience the over-smoothing problem, because the model converges to a meaningful representation as $t \to \infty$ given by a fixed point \cite{pmlr-v119-xhonneux20a, chamberlain2021grand, poli2019graph}.

The Continuous Graph Neural Network (CGNN) \cite{pmlr-v119-xhonneux20a}  architecture consists of three main components. First, a neural encoder projects each node into a latent space based on its features, resulting in $H_0 = \mathcal{E}(X)$, where $\mathcal{E}$ represents the encoder. Next, $H_0$ is used as the initial value $H(0)$, and an ODE is designed to define the continuous dynamics of node representations, effectively modeling long-term dependencies between nodes. Finally, the node representations at the end time $t_n$ (i.e., $H(t_n)$) are used for downstream tasks through a decoder $\mathcal{D}$, yielding a node-label matrix $Y = \mathcal{D}(H(t_n))$.

One example of an ODE to model the continuous dynamics of node representations is 
\begin{equation}
    \label{eq:CGNN_ODE}
    \dot Q(t) = (\Tilde{A}- \mathds{1}) Q(t) + Q(t)( W - \mathds{1}) + H_0,
\end{equation}
where $W \in \mathbb{R}^{ d\, \times \, d}$ is a weight matrix. This is a differential Sylvester equation where the left and right operators govern system dynamics with the inhomogeneity being the initial input $E$ from the encoder network. This system includes both extensions mentioned earlier: non-commuting operators and a non-square system. We integrate this system using our matrix ETD solver and show that it can solve complex high-dimensional Sylvester systems.

\textbf{Datasets.}   We study our schemes using the three most popular citation networks in the literature, Cora \cite{mccallum2000automating}, CiteSeer \cite{sen2008collective} and Pubmed \cite{namata2012query}. We make use of random weight initializations of fixed splits of these datasets and show a direct comparison in Table \ref{tab:GNN_ETD}.

\textbf{Implementation details.} All schemes in this example have been implemented in PyTorch \cite{paszke2019pytorch} using torchdiffeq \cite{chen2018neural} and PyTorch geometric \cite{fey2019fast}. The hyperparameters from Xhonneux et al. \cite{pmlr-v119-xhonneux20a} were used for all schemes. We do, however, use Theorem \ref{thm:BCH_convergence} to make sure our BCH series is convergent by varying the hyperparameter $\alpha$. To ensure non-singularity of $Z_h$ in our schemes during training, we add some regularization to $Z_h$ such that we do not have any null eigenvalues. We also use the ANODE augmentation scheme \cite{anode} to stabilize training.

To demonstrate that our METD schemes work on differential Sylvester equations, we solve the ODE \eqref{eq:CGNN_ODE} and show the test accuracy results for the METD schemes compared to the default in most packages which is the adaptive Runge-Kutta 4(5) Dormand-Price (Dopri5) scheme which we take as the baseline. We also compare our schemes to the explicit 4th order Adams–Bashforth (Explicit AB) and the implicit 4th order Adams–Moulton (Implicit AM) schemes. Table \ref{tab:GNN_ETD} summarises that we achieve comparable results to the original CGNN implementation using the Dopri5 solver which shows that the extended METD schemes work well with differential Sylvester equations. 

\begin{table}[t]
\caption{Node classification test accuracy and std for 20 random initializations using the Planetoid train-val-test splits. Bold results indicate the highest test accuracy for its respective dataset. \textbf{$^*$}Results obtained by running authors' original code with parameters given in their paper.}
\vskip 0.15in
\begin{center}
\begin{small}
\begin{sc}
\begin{tabular}{lcccr}
\hline
  Planetoid splits & \bf CORA & \bf CiteSeer & \bf PubMed \\
  $\#$Nodes & 2,708 &  3,327 &  18,717 \\
  $\#$Edges & 5,278 &  4,676 &  44,327 \\ 
  $\#$Classes & 7 &  6 &  3 \\[0.5ex]
\hline
\bf Dopri5$^*$ & 81.8 $\pm$ \,0.7 & 68.1 $\pm$ 1.2 & 80.3 $\pm$ 0.3 \\ 
\bf Explicit AB & 32.0 $\pm$ \,1.9 & 30.1 $\pm$ 1.8 & 31.1 $\pm$ 1.9 \\
\bf Implicit AM & 82.3 $\pm$ \,0.6 & 69.0 $\pm$ 0.7 & 80.5 $\pm$ 0.4 \\
\hline
 \bf METD1$\_$BCH1 & 80.4 $\pm$ \,1.4 & 67.3 $\pm$ 1.0 & 80.0 $\pm$ 0.5 \\
 \bf METD1$\_$BCH2 & 82.8 $\pm$ \,0.5 & 67.5 $\pm$ 1.2 & 80.4 $\pm$ 0.7 \\
 \bf METD1$\_$BCH3 & \textbf{83.0 $\pm$ 0.5} & 68.9 $\pm$ 1.1 & \textbf{80.9 $\pm$ 0.6} \\
 \bf METD2$\_$BCH1 & 80.8 $\pm$ \,1.4 & 67.9 $\pm$ 1.1 & 80.3 $\pm$ 0.8\\
 \bf METD2$\_$BCH2 & 83.0 $\pm$ \,0.6 & 68.0 $\pm$ 1.0 & 80.5 $\pm$ 0.6\\
 \bf METD2$\_$BCH3 & 82.9 $\pm$ \,0.5 & \textbf{69.1 $\pm$ 0.8} & 80.8 $\pm$ 0.5\\
\hline
\end{tabular}
\end{sc}
\end{small}
\end{center}
\vskip -0.1in
\label{tab:GNN_ETD}
\end{table}

\begin{table}[t]
\caption{Average running times per epoch for various fixed time step schemes on the Cora dataset. Out of the step sizes tried for the explicit AB scheme, only $\Delta t= 0.005$ was stable. The implicit AM scheme is the most expensive in terms of runtime per epoch compared to other schemes. The METD1$\_$BCH3 scheme has the lowest running time.}
\vskip 0.15in
\begin{center}
\begin{small}
\begin{sc}
\begin{tabular}{lcccr}
\hline
  Fixed time step schemes & \bf Explicit AB & \bf Implicit AM & \bf METD1$\_$BCH3 \\
  & ($\Delta t = 0.005$) & ($\Delta t = 0.01$) &  ($\Delta t = 0.01$) \\
  \hline
\bf Average time per epoch (s) & 9.2 & 19.2 &  2.9 \\ 

\hline
\end{tabular}
\end{sc}
\end{small}
\end{center}
\vskip -0.1in
\label{tab:runnung_times}
\end{table}

The explicit AB scheme is unstable for all step sizes except for very small time steps ($\Delta t = 0.005$) which makes using this solver impractical since a single epoch takes a significant amount of time. The implicit AM scheme is stable for all step sizes except large ones ($\Delta t = 10$) however, with a large step size, the implicit equations can become more difficult to solve. Table \ref{tab:runnung_times} shows a comparison of running times per epoch for some schemes. We see that the METD1$\_$BCH3 scheme is faster as well as just as accurate (in terms of test accuracy from \cref{tab:GNN_ETD}) as the implicit AM scheme.

Note that we are not aiming to outperform the benchmark, but rather to show that our schemes can integrate these complex systems and perform just as well as state-of-the-art integrators while taking less time.

In some cases, we do actually outperform the original implementation. For example, the METD1$\_$BCH3 scheme outperforms the original implementation for all the benchmark datasets tested. We also observed that using METD, we can choose the time step to be significantly larger than when using other methods. With just $3$ layers, we can achieve comparable accuracy to the Dopri5 CGNN, which has potentially hundreds of layers depending on stability. As noted in Remark~\ref{rmk:order_limit}, due to the linear freezing approximation, both METD1\_BCH$n$ and METD2\_BCH$n$ achieve at most second-order accuracy regardless of $n$. However, higher BCH truncation depths can reduce the error constants, as reflected in the improved test accuracies observed for larger $n$ in Table~\ref{tab:GNN_ETD}.

\section{Conclusion}

We have extended the ETD numerical method to matrix differential equations which have a stiff linear part using elements from Lie algebra. A significant advantage of these schemes stems from the fact that the linear part is solved exactly, which makes these stiff problems tractable from a numerical point of view. Integrating matrix-valued dynamical systems using vectorization and classic ETD is prohibitively expensive for large systems, while METD schemes translate the ideas of ETD onto matrix evolution equations without a large increase in computational effort.

In terms of algorithms, we derived first and second-order METD multistep schemes as well as a Runge–Kutta form of second-order which does not require an initialization. These METD schemes enjoy desirable stability properties and can be applied to a variety of problems including matrix-ODE and PDE systems. In addition, we introduced an explicit order-$p$ METD$p$ family based on backward-interpolation of the nonlinear term and a controlled truncation of nested commutator corrections. We proved that, under standard semigroup stability bounds for the linear propagators and regularity assumptions on the nonlinearity, the resulting schemes are consistent with local defect $\mathcal{O}(h^{p+1})$ and achieve global convergence of order $p$ on finite time intervals. These results provide a rigorous foundation for using higher-order METD methods in stiff matrix evolution equations without resorting to vectorization.

We have carried out a number of numerical tests for the presented METD schemes. We validated the predicted convergence orders on differential Lyapunov systems and on a stiff nonlinear Allen–Cahn benchmark obtained from PDE discretization. On the Allen–Cahn example, METD decouples stiffness from the time step by integrating the Laplacian exactly, enabling a practical accuracy–cost tradeoff that is difficult for standard explicit adaptive Runge–Kutta methods on this regime. In these examples, we only consider the case of commuting left and right operators, which is typical when arriving at the equations from stochastic systems. We also showed how these methods succeed in the integration of a severely stiff system that models jet formation in atmospheres of large planets which would usually require an extremely small time step.

Finally, we considered a continuous GNN where a Sylvester equation replaces the discrete layers typically seen in neural networks. In order to apply the developed methods for this problem, it was necessary to extend the METD scheme to non-commuting operators with a non-square system and to relax some of the assumptions about the involved commutators. We reported test accuracies for three benchmark datasets and compared these to the default adaptive Runge-Kutta 4(5) Dormand-Price scheme as well as other fixed time step schemes showing comparable accuracy. We show a significant speedup compared to other fixed time step schemes. In our tests we observed that METD permits larger step sizes (fewer solver steps) at comparable accuracy.

\appendix
\section{Useful lemmas and derivations}

\subsection{Useful lemmas}

\begin{lem}
\label{lma1}
Let $A$ be an arbitrary matrix $A \in \mathbb{C}^{n \times n}$ and $t \in \mathbb{R}$, then the commutator $[A,e^{-tA}] = 0$.
\end{lem}

\begin{pf}

 We have that $[A,e^{-tA}] = Ae^{-tA} - e^{-tA}A \,$ and expanding both matrix exponentials as a power series gives us 
\begin{align*}
     Ae^{-tA} - e^{-tA}A &= A \left[\sum_{n=0}^{\infty} \frac{(-t)^n}{n!}(A)^{n} \right] - \left[ \sum_{n=0}^{\infty} \frac{(-t)^n}{n!}(A)^{n}\right] A\\ 
     &= \sum_{n=0}^{\infty} \frac{(-t)^n}{n!}(A)^{n+1} - \sum_{n=0}^{\infty} \frac{(-t)^n}{n!}(A)^{n+1}\\
     &= 0.
\end{align*}
\end{pf}

\begin{lem}
\label{order_2_METD1}
Let $L \in \mathbb{C}^{n \times n}$ then $\int_0^{t} s e^{sL}\,ds = \mathcal O( t^2).$ 
\end{lem}
\begin{pf}
By integrating, we have
\begin{align*}
        \int_0^{t} s e^{sL}\,ds &= L^{-2}\left( te^{tL}L - e^{tL} + \mathds{1}  \right)\\
        &= L^{-2}\left(t \sum_{n=0}^\infty \frac{t^n}{n!}L^{n+1} - \sum_{n=0}^\infty \frac{t^n}{n!}L^n + \mathds{1} \right)\\
        &=L^{-2}  \left( [tL + t^2L^{2} + \frac{t^3}{2}L^3 + \dots ] - 
        [\mathds{1} + tL + \frac{t^2}{2}L^2 + \dots] + \mathds{1}      \right)
        \\&= L^{-2} \left( \frac{t^2}{2}L^2 +  \mathcal O( t^3)     \right)\\
        &= \mathcal O( t^2)
\end{align*}
\end{pf} 

\begin{lem}
\label{lem:ad_R}
For $R,X\in\RR^{n\times n}$ and $s\in\RR$,
\begin{equation*}
    e^{s\,\ad_R}(X)=e^{-sR}Xe^{sR},
\qquad \ad_R(X):=[X,R].
\end{equation*}
\end{lem}
\begin{pf}
    Define $Y(s) = e^{-sR}Xe^{sR}$. Differentiate
    \begin{align*}
        Y'(s) &= (-R)e^{-sR}Xe^{sR} + e^{-sR}X(Re^{sR})\\
        &=-RY(s) + Y(s)R\\
        &= Y(s)R - RY(s)\\
        &= [Y(s),R] \\
        &= \ad_R(Y(s))
    \end{align*}
    so $Y$ solves the linear ODE
    \begin{equation}
        Y'(s) = \ad_R(Y(s)), \quad Y(0) = X.
    \end{equation}
    The solution of $Y' = \ad_R(Y)$ with $Y(0) = X$ is exactly
    \begin{equation}
        Y(s) = e^{s\,\ad_R}(X),
    \end{equation}
    hence 
    \begin{equation}
        e^{s\,\ad_R}(X) = e^{-sR}Xe^{sR}.
    \end{equation}
\end{pf}

\begin{lem}
\label{lem:diff-deriv}
Let $(\mathcal{B},\|\cdot\|)$ be a Banach space and let $G\in C^{m}([t_n-mh,t_n];\mathcal{B})$.
Define the backward difference with step size $h>0$ by
\begin{equation*}
(\nabla_h G)(t):=G(t)-G(t-h),\qquad \nabla_h^{0}G:=G,\qquad \nabla_h^{m}:=\nabla_h(\nabla_h^{m-1}).
\end{equation*}
Then $\nabla_h^{m}G(t_n)$ admits the integral representation
\begin{equation}
\label{eq:diff-int-rep}
\nabla_h^{m}G(t_n)
= \int_{[0,h]^m} G^{(m)}\!\left(t_n-\sum_{i=1}^m s_i\right)\,ds_1\cdots ds_m,
\end{equation}
and consequently the norm bound
\begin{equation}
\label{eq:diff-deriv-bound}
\|\nabla_h^{m}G(t_n)\|
\le h^{m}\,\sup_{t\in[t_n-mh,t_n]}\|G^{(m)}(t)\|.
\end{equation}
In particular, if $\sup_{t\in[t_n-mh,t_n]}\|G^{(m)}(t)\|$ is bounded independently of $h$,
then $\nabla_h^{m}G(t_n)=\mathcal{O}(h^{m})$ as $h\to 0$.
\end{lem}

\begin{pf}
For $m=1$, the fundamental theorem of calculus gives
\begin{equation*}
\nabla_h G(t)=G(t)-G(t-h)=\int_{t-h}^{t} G'(s)\,ds=\int_{0}^{h} G'(t-s)\,ds.
\end{equation*}
Assume \eqref{eq:diff-int-rep} holds for $m-1$. Apply the $m=1$ identity to the function
$F:=\nabla_h^{m-1}G$
\begin{equation}
\label{eq:diff-integ}
\nabla_h^{m}G(t)=\nabla_h F(t)=\int_{0}^{h} F'(t-s_m)\,ds_m.
\end{equation}
Since differentiation commutes with $\nabla_h$ on $C^{m}$ functions (linearity plus the chain rule),
we have $F'=\nabla_h^{m-1}G'$. Using the induction hypothesis with $G'$ in place of $G$ yields
\begin{equation}
\label{eq:m-integral}
F'(t-s_m)
= \int_{[0,h]^{m-1}} G^{(m)}\!\left(t-s_m-\sum_{i=1}^{m-1}s_i\right)\,ds_1\cdots ds_{m-1}.
\end{equation}
Substituting \eqref{eq:m-integral} into \eqref{eq:diff-integ} yields
\begin{align*}
\nabla_h^{m}G(t)
&=\int_{0}^{h}\left(\int_{[0,h]^{m-1}}
G^{(m)}\!\left(t-s_m-\sum_{i=1}^{m-1}s_i\right)\,ds_1\cdots ds_{m-1}\right)ds_m\\
&=\int_{[0,h]^m} G^{(m)}\!\left(t-\sum_{i=1}^{m}s_i\right)\,ds_1\cdots ds_m,
\end{align*}
which is \eqref{eq:diff-int-rep}.

Taking norms on both sides gives
\begin{equation*}
\|\nabla_h^{m}G(t_n)\|
\le \int_{[0,h]^m}\left\|G^{(m)}\!\left(t_n-\sum_{i=1}^{m}s_i\right)\right\|\,ds_1\cdots ds_m.
\end{equation*}
For all $(s_1,\dots,s_m)\in[0,h]^m$ we have
$t_n-\sum_{i=1}^{m}s_i \in [t_n-mh,t_n]$, hence
\begin{equation*}
\left\|G^{(m)}\!\left(t_n-\sum_{i=1}^{m}s_i\right)\right\|
\le \sup_{t\in[t_n-mh,t_n]}\|G^{(m)}(t)\|.
\end{equation*}
Therefore
\begin{equation*}
\|\nabla_h^{m}G(t_n)\|
\le \left(\sup_{t\in[t_n-mh,t_n]}\|G^{(m)}(t)\|\right)
\int_{[0,h]^m} 1\,ds_1\cdots ds_m.
\end{equation*}
Finally, the integral of $1$ over the $m$-dimensional box $[0,h]^m$ equals its volume,
which is $h^m$. This proves \eqref{eq:diff-deriv-bound}.
\end{pf}

\begin{lem}
\label{matrix_exp_identity}
Let $A \in \mathbb{C}^{n \times n}$ and $\Tilde{A} \in \mathbb{C}^{m \times m}$ where $m>n$ with $A$ zero-padded to form $\Tilde{A}$ i.e. $\Tilde{A} = \begin{bmatrix}
  A & 0_{n \times (m-n)}\\ 
  0_{(m-n)\times n} & 0_{(m-n)}
\end{bmatrix}$. The matrix exponential of $\Tilde{A}$ will then be $e^{\Tilde{A}} = \begin{bmatrix}
  e^{A} & 0\\ 
  0 & \mathds 1_{(m-n)}
\end{bmatrix}$, i.e. $e^{A}$ padded by an identity matrix $\mathds 1_{(m-n)} $.
\end{lem}
\begin{pf} 
Using the definition of the matrix exponential, 
\begin{align*}
   e^{\Tilde{A}} &=  \sum_{k=0}^{\infty} \frac{1}{k!}(\Tilde{A})^{k} 
   \\ &= \begin{bmatrix}
  \mathds 1_{n} & 0\\ 
  0 & \mathds 1_{(m-n)}\end{bmatrix} + \begin{bmatrix}
  A & 0_{n \times (m-n)}\\ 
  0_{(m-n)\times n} & 0_{(m-n)}\end{bmatrix} + \frac{1}{2!}\begin{bmatrix}
  A^2 & 0_{n \times (m-n)}\\ 
  0_{(m-n)\times n} & 0_{(m-n)}\end{bmatrix} +\dots 
  \\ &= \begin{bmatrix}
  \mathds 1_{n} + A +  \frac{1}{2!}A^2 +\dots & 0\\ 
  0 & \mathds 1_{(m-n)}\end{bmatrix}
  \\ &= \begin{bmatrix}
  e^{A} & 0\\ 
  0 & \mathds 1_{(m-n)}
  \end{bmatrix}.
\end{align*}

\end{pf} 

\subsection{Derivation of METD1}
\label{ETD1_deriv}
Take the first-order approximation of the integral of the nonlinear term
\begin{equation*}
    \int_0^{\Delta t} e^{(\Delta t-\tau)(L+R)} \,d\tau N(t_n),
\end{equation*}
and set $A=: L+R$ with a substitution $s = \Delta t-\tau$. Changing the limits of integration with the substitution gives us the integral with solution 
\begin{align*}
    \int_0^{\Delta t} e^{sA} \,ds \,N(t_n) &= [\left( e^{\Delta t \,A} - \mathds{1}\right)A^{-1}]N(t_n)
    \\ &= [(L+R)^{-1} \left( e^{\Delta t \,(L+R)} - \mathds{1}\right)]N(t_n)
\end{align*}
with $\left( e^{\Delta t \,A} - \mathds{1}\right)$ and $A^{-1}$ commuting from Lemma~\ref{lma1}.

\subsection{Derivation of METD2}
\label{app:METD_deriv}
By setting $\xi_n = (N_{n}- N_{n-1})/\Delta t$, we can arrive at the following integrals 
\begin{align*}
    \int_0^{\Delta t} e^{(\Delta t-\tau)L} N(t_n+\tau) e^{(\Delta t-\tau)R}\,d\tau =&\int_0^{\Delta t} e^{(\Delta t-\tau)(L+R)} N(t_n+\tau)\,d\tau \\ &+ \int_0^{\Delta t} (\Delta t-\tau) e^{(\Delta t-\tau)L} [N(t_n+\tau),R]\,d\tau + \mathcal O(\Delta t^3) \\= 
    &\int_0^{\Delta t} e^{(\Delta t-\tau)(L+R)} (N_n + \tau \xi_n)\,d\tau \, \\ &+ \int_0^{\Delta t} (\Delta t-\tau) e^{(\Delta t-\tau)L} [(N_n + \tau \xi_n),R]\,d\tau \\=& I_n + C_n
\end{align*}
where we are solving both integrals separately for simplicity.
\begin{align*}
    I_n &= \int_0^{\Delta t} e^{(\Delta t-\tau)(L+R)}\,d\tau\, N_n + \int_0^{\Delta t} e^{(\Delta t-\tau)(L+R)} \tau \,d\tau\, \xi_n \\
    &= (L+R)^{-1}\left(e^{\Delta t (L+R)} - \mathds{1}\right)\,N_n + (L+R)^{-2}\left(e^{\Delta t (L+R)} - \Delta t (L+R) - \mathds{1}\right)\xi_n
\end{align*}
and
\begin{align*}
    C_n &=  \int_0^{\Delta t} (\Delta t-\tau) e^{(\Delta t-\tau)L} [(N_n + \tau \xi_n),R]\,d\tau\\
    &= \int_0^{\Delta t} (\Delta t-\tau) e^{(\Delta t-\tau)L}\,d\tau\, [N_n,R] \, + \int_0^{\Delta t} \tau (\Delta t-\tau) e^{(\Delta t-\tau)L}\,d\tau\, [\xi_n,R]\\
    &= L^{-2} \left( \Delta t e^{\Delta t L}L -  e^{\Delta t L} + \mathds{1}  \right)  [N_n,R] \,+
    L^{-3} \left(  \Delta t L + \Delta t e^{\Delta t L}L -  2e^{\Delta t L} + 2 \mathds{1}  \right)[\xi_n,R].
\end{align*}
Note that the $[\xi_n,R]$ term is $O(\Delta t^3)$ and should be dropped for METD2.

\subsection{Derivation of METD$p$}
\label{app:METDp_deriv}
Let us restate the system that we want to solve which is
\begin{equation}
    \dot Q(t) = LQ + QR + N(Q,t)
\end{equation}
with $Q,L,R \in \RR^{n \times n}$. We assume $[L,R] = LR - RL =0$ and lets set $h=\Delta t$. One step, $t_{n+1} = t_n + h$, we then have the variation of constant formula
\begin{equation}
\label{eq:var_constants}
    Q_{n+1} = e^{hL}Q_ne^{hR} + \underbrace{\int_0^h e^{(h-\tau)L}\, N(t_n + \tau)e^{(h-\tau)R}\, d\tau}_\text{:=I} .
\end{equation}
For any matrix $R$, we have a linear operator on the space of matrices 
\begin{equation}
    \ad_R : \RR^{n \times n} \to \RR^{n \times n}, \qquad \ad_R(X) := [X,R] = XR - RX.
\end{equation}
We can have that 
\begin{align*}
    \ad_R^{\,0}(X) &= X,\\
    \ad_R^{\,1}(X) &= [X,R],\\
    \ad_R^{\,2}(X) &= \ad_R(\ad_R(X)) = [[X,R],R],\\
    \ad_R^{\,3}(X) &= [[[X,R],R],R],\\
\end{align*}
so $\ad_R^{\,j}(X)$ is the $j^{th}$-times nested commutator with R. We can expand $\ad_R$ in a power series
\begin{equation}
    e^{s\,\ad_R} = \sum_{j=0}^{\infty}\frac{s^j}{j!}\ad_R^{\,j}
\end{equation}
so 
\begin{align}
    e^{s\,\ad_R}(X) = \sum_{j=0}^{\infty}\frac{s^j}{j!}\ad_R^{\,j} (X) &= \ad_R^{\,0}(X) + s\, \ad_R^{\,1}(X)+ \frac{s^2}{2}\,\ad_R^{\,2}(X) + ...\\
    &= X + s[X,R] + \frac{s^2}{2}[[X,R],R] + ...
\end{align}
so $e^{s\,\ad_R}(X)$ is a matrix obtained by adding nested commutators.\\
Now lets us set $s:= h-\tau \geq 0$, we can manipulate the term in the integral in I in \eqref{eq:var_constants} like so 
\begin{align}
    e^{sL}Ne^{sR} &= e^{sL}(e^{sR}e^{-sR})Ne^{sR}\\
    &=e^{sL}e^{sR}[e^{-sR}Ne^{sR}]
\end{align}
since we have that $[L,R] =0$, we have 
\begin{equation}
\label{eq:1_no_ad}
    e^{sL}Ne^{sR} = e^{s(L+R)}(e^{-sR}Ne^{sR}).
\end{equation}

Using Lemma~\ref{lem:ad_R}, we can write \eqref{eq:1_no_ad} as
\begin{equation}
\label{eq:2_ad}
    e^{sL}Ne^{sR} = e^{s(L+R)}(e^{-sR}Ne^{sR}) = e^{s(L+R)}e^{s\,\ad_R}(N).
\end{equation}
We can now substitute \eqref{eq:2_ad} into the integral term in \eqref{eq:var_constants} so that we have 
\begin{align}
    \textit{I} &= \int_0^h e^{(h-\tau)L}\, N(t_n + \tau)e^{(h-\tau)R}\, d\tau \\
    &=\int_0^h e^{(h-\tau)(L+R)}\, e^{(h-\tau)\,\ad_R}\, (N(t_n + \tau))\, d\tau.
    \label{eq:1_int}
\end{align}

Now we need to do a two things:

\begin{itemize}
    \item We need to approximate the nonlinear term, $N(t_n + \tau)$, by a polynomial of degree $(p-1)$ in $\tau$.
    \item We need to expand $e^{(h-\tau)\,\ad_R}\,(\cdot)$ in a series and truncate so the local error is $O(h^{p+1})$.
\end{itemize}
Following Cox et al.\cite{cox-matthews:2002} for the nonlinear term, we have the backward difference operator $\nabla$ which acts on the discrete sequence $N_n := N(Q_n,t_n)$ where
\begin{equation}
    \nabla^0N_n = N_n, \quad \nabla N_n = N_n- N_{n-1}, \quad \nabla^2 N_n = N_n - 2N_{n-1} + N_{n-2}
\end{equation}
so 
\begin{equation}
    \nabla^mN_n = \sum_{\ell=0}^{m}(-1)^{\ell} \binom{m}{\ell} N_{n-\ell}.
\end{equation}
Note that throughout this derivation, we use $\nabla := \nabla_h$ for the backward difference with step size $h$.
The degree $(p-1)$ interpolant of the function $t \mapsto N(t)$ at grid points $t_n,t_{n-1}, ..., t_{n-p+1}$ evaluated at $t_n+\tau$ is then given by
\begin{equation}
\label{eq:nonlinear_approx}
    N(t_n+\tau) = \sum_{m=0}^{p-1} (-1)^m \binom{-\theta}{m}\nabla^mN_n + O(h^p), \quad 0 \leq \tau \leq h,
\end{equation}
where $\theta h =\tau,\, \theta \in [0,1]$. We now write the integral \eqref{eq:1_int} like
\begin{equation}
\label{eq:int_2}
    \textit{I} = h\int_0^1 e^{h(1-\theta)(L+R)}\, e^{h(1-\theta)\,\ad_R}\, (N(t_n + h\theta))\, d\theta \,+ O(h^{p+1})\,,
\end{equation}
where we have done a change of variables $\tau = \theta h$ so $d\tau = h\, d\theta$ and $h-\tau = h(1-\theta)$. We can set $A:=h(L+R)$ and insert \eqref{eq:nonlinear_approx} into \eqref{eq:int_2} to get
\begin{equation}
\label{eq:int_3}
        \textit{I} = h\int_0^1 e^{(1-\theta)A}\, \underbrace{e^{h(1-\theta)\,\ad_R}\, \left(\sum_{m=0}^{p-1} (-1)^m\binom{-\theta}{m}\nabla^mN_n\right)}_\text{:=II}\, d\theta \, \,+ O(h^{p+1})\, .
\end{equation}

Now we need to expand II in \eqref{eq:int_3} and truncate. We have that
\begin{equation}
\label{eq:expansion-exp-adr}
    e^{h(1-\theta)\,\ad_R}\, (X) = \sum_{j=0}^{\infty} \frac{h^j(1-\theta)^j}{j!}\,\ad_R^{\,j}\,(X).
\end{equation}
If we substitute $X = \binom{-\theta}{m}\nabla^mN_n$ into \eqref{eq:expansion-exp-adr}, we have
\begin{equation}
\label{eq:expansion_exp_adr_N}
    e^{h(1-\theta)\,\ad_R}\, \left((-1)^m\binom{-\theta}{m}\nabla^mN_n \right) = (-1)^m\binom{-\theta}{m}\sum_{j=0}^{\infty} \frac{h^j(1-\theta)^j}{j!}\,\ad_R^{\,j}\,\left(\nabla^mN_n \right).
\end{equation}
So now if we substitute \eqref{eq:expansion_exp_adr_N} into \eqref{eq:int_3}, we get
\begin{equation}
    \textit{I} = h\sum_{m=0}^{p-1}\sum_{j=0}^{\infty}\frac{(-1)^mh^j}{j!}\left[\int_0^1 e^{(1-\theta)A}(1-\theta)^j\binom{-\theta}{m}\,d\theta\right]\,\ad_R^{\,j}\,\left(\nabla^mN_n \right)\, + O(h^{p+1})\,.
\end{equation}

By Lemma~\ref{lem:diff-deriv}, if
$t\mapsto N(Q(t),t)\in C^{p}$ along the exact solution, then
$\|\nabla^m N_n\|=\mathcal{O}(h^m)$. 
So we have contributions from $h$, $h^j$ and $\nabla^mN_n$ so $h \cdot h^j \cdot \nabla^mN_n \sim h^{1+j+m}$. Then to get local truncation error of $O(h^{p+1})$, we keep all terms with
\begin{equation}
    m+j+1 \leq p \quad \Leftrightarrow \quad m+j \leq p-1 \,,
\end{equation}
and drop the rest so 
\begin{equation}
\label{eq:I_p}
    \textit{I}_p = h\sum_{m=0}^{p-1}\sum_{j=0}^{p-1-m}h^j\,C_{m,j}(A)\,\ad_R^{\,j}\,\left(\nabla^mN_n \right)\,,
\end{equation}
where we have the integral coefficients defined as
\begin{equation}
\label{eq:c_mj}
    C_{m,j}(A):=\frac{(-1)^m}{j!}\int_0^1 e^{(1-\theta)A}(1-\theta)^j\binom{-\theta}{m}\,d\theta\,.
\end{equation}
Now we need to reduce these coefficients $C_{m,j}(A)$ to linear combinations of $\varphi_k(A)$ where
\begin{equation}
    \varphi_k(A) = \int_0^{1} e^{(1 - \theta)A} \frac{\theta^{k-1}}{(k-1)!}\, d\theta\,.
\end{equation}

For fixed $m,j$, the scalar polynomial $(1-\theta)^j\binom{-\theta}{m}$ is a polynomial in $\theta$ of degree $m+j$. So expanding the scalar polynomial gives
\begin{equation}
\label{eq:theta_expansion}
    (1-\theta)^j\binom{-\theta}{m} = \sum_{q=0}^{m+j}\alpha_{m,j,q}\, \theta^q \,,
\end{equation}
where $\alpha_{m,j,q}$ is the coefficient of $\theta^q$. We also have that 
\begin{equation}
\label{eq:var_phi_q}
    \int_0^1e^{(1-\theta)A}\theta^q\,d\theta = q!\varphi_{q+1}(A).
\end{equation}
So substituting \eqref{eq:theta_expansion} and \eqref{eq:var_phi_q} into \eqref{eq:c_mj} to get 
\begin{equation}
\label{eq:c_mj_theta}
    C_{m,j}(A) = \frac{(-1)^m}{j!}\sum_{q=0}^{m+j} \alpha_{m,j,q}\,q!\,\varphi_{q+1}(A).
\end{equation}
Now if we combine \eqref{eq:c_mj_theta}, \eqref{eq:I_p} and \eqref{eq:var_constants}, we get the final arbitrary-order METD$p$ scheme

\begin{equation}
    \label{eq:final_schema}
    Q_{n+1} = e^{hL}Q_ne^{hR} +
    h\sum_{m=0}^{p-1}\sum_{j=0}^{p-1-m}h^j\,\left( \frac{(-1)^m}{j!}\sum_{q=0}^{m+j} \alpha_{m,j,q}\,q!\,\varphi_{q+1}(A) \right)\,\ad_R^{\,j}\,\left(\nabla^mN_n \right)\,.
\end{equation}
We now propose algorithm~\ref{alg:build_cmj_algorithmic}  for the construction of the coefficient matrices $C_{m,j}$ in the METD$p$ update.  

\begin{algorithm}[t]
\caption{Build coefficient matrices $C_{m,j}(A)$}
\label{alg:build_cmj_algorithmic}
\begin{algorithmic}[1]
\REQUIRE Order $p$, $\varphi_k(A)$ for $k=1,\ldots,p$
\ENSURE $C_{m,j}(A)$ for $0\le m\le p-1$, $0\le j\le p-1-m$
\STATE Precompute $f_r=r!$ for $r=0,\ldots,p$
\STATE Build and store polynomial coefficients $b^{(m)}$ of $\binom{-\theta}{m}$ for $m=0,\ldots,p-1$
\FOR{$m=0,\ldots,p-1$}
    \FOR{$j=0,\ldots,p-1-m$}
        \STATE Build $u^{(j)}$ for $(1-\theta)^j=\sum_{k=0}^{j}u_k^{(j)}\theta^k$, with $u_k^{(j)}=\binom{j}{k}(-1)^k$
        \STATE $\alpha \leftarrow \mathrm{conv}(u^{(j)},b^{(m)})$
        \STATE $S \leftarrow 0$
        \FOR{$q=0,\ldots,m+j$}
            \STATE $S \leftarrow S + \alpha_q\,f_q\,\varphi_{q+1}(A)$
        \ENDFOR
        \STATE $C_{m,j}(A) \leftarrow \dfrac{(-1)^m}{f_j}\,S$
    \ENDFOR
\ENDFOR
\RETURN $\{C_{m,j}(A)\}$
\end{algorithmic}
\end{algorithm}

\bibliographystyle{elsarticle-num} 
\bibliography{bibliography}

\end{document}